\documentclass[journal,final]{IEEEtran}

\ifCLASSINFOpdf
\else
\fi

\usepackage[cmex10]{amsmath}
\usepackage{algorithmic}
\usepackage{bm, cite, psfrag, graphics, epsfig, latexsym, amsmath, color, amsfonts, amssymb, graphicx, algorithm, algorithmic,mathrsfs}
\usepackage{xcolor}
\usepackage{multirow}
\usepackage{wasysym}
\usepackage{pgf}
\usepackage{subcaption}\usepackage{array}
\usepackage{url}

\newtheorem{theorem}{Theorem}[section]
\newtheorem{lemma}{Lemma}[section]
\newtheorem{definition}{Definition}[section]
\newtheorem{remark}{Remark}[section]

\begin{document}

\title{Iterative learning and extremum seeking for repetitive time-varying mappings}

\author{Zhixing Cao, Hans-Bernd D{\"u}rr, Christian Ebenbauer, Frank Allg{\"o}wer, Furong Gao$^*$
\thanks{Zhixing Cao and Furong Gao are with the Department of Chemical and Biomolecular Engineering, Hong Kong University of Science and Technology, Hong Kong SAR. Zhixing Cao is also with Harvard John A. Paulson School of Engineering and Applied Sciences, Harvard University, Cambridge, MA 02138, United States. Hans-Bernd D{\"u}rr, Christian Ebenbauer and Frank Allg{\"o}wer are with the Institute for Systems Theory and Automatic Control, University of Stuttgart, 70569 Germany. Furong Gao is also with Fok Ying Tung Research Institute, Hong Kong University of Science and Technology, Guangzhou, China. E-mail: edwardcao@g.harvard.edu, \{hans-bernd.duerr, ce, frank.allgower\}@ist.uni-stuttgart.de, kefgao@ust.hk.

This work is supported by the National Natural Science Foundation of China $\#$61227005, Guangdong Innovative and Entrepreneurial Research Team Program $\#$2013G076, Guangzhou Science and Technology Bureau Project $\#$12190007.
}
}

\markboth{IEEE TRANSACTIONS ON AUTOMATIC CONTROL}%
{Cao \MakeLowercase{\textit{et al.}}: Iterative learning extremum seeking for unknown reference tracking via Lie bracket approximation}



\maketitle

\begin{abstract}
In this paper, we develop an {{iterative learning control}} method integrated with {{extremum seeking control}} to track a time-varying optimizer within {finite time horizon}. The behavior of the extremum seeking system is analyzed via an approximating system -- the modified Lie bracket system. The modified Lie bracket system is essentially an online integral-type iterative learning control law. The paper contributes to two fields, namely, iterative learning control and extremum seeking. First, an online integral type iterative learning control with a forgetting factor is proposed. Its convergence is analyzed via $k$-dependent (iteration-dependent) contraction mapping in a Banach space equipped with {{so called}} $\lambda$-norm. Second, the iterative learning extremum seeking system can be {{interpreted as an iterative learning control with the approximation error as ``disturbance''.}} The tracking error of its modified Lie bracket system can be shown uniformly bounded in terms of iterations by selecting a sufficiently large dither frequency. Furthermore, it is shown that the tracking error will eventually converge to a set.  {{The center of the set corresponds to the limit solution of the ``disturbance-free'' system, and its radius can be controlled by the frequency.}}%The convergence of the modified Lie bracket system is concluded by contraction mapping in a Banach space equipped with $\lambda$-norm. Moreover, it can be shown that the original system is bounded in the neighborhood of the time varying optimizer uniformly in terms of iterations.
\end{abstract}


% \fi
\IEEEpeerreviewmaketitle

\section{Introduction}

Considerable research efforts have been devoted to extremum seeking (ES) control {{over}} the last several decades. The mechanism of extremum seeking control is to optimize a certain system performance measure (cost function) by adaptively adjusting the system parameters merely according to output measurements of the plant. Since there is little knowledge required about the plant dynamics, extremum seeking control has attracted the attention from various engineering domains, e.g., bioreactor, combustion, compressor \cite{Krstic_book,Tan2010,Dochain}. 

{Both the success of extremum seeking in industrial application and the uniqueness of repetitive process optimization (RPO) problem (e.g., polymer-melt-front-velocity (PMFV) optimization in Section \ref{PMFV}) motivated this study. Theoretically, the RPO problem has many different features compared to the standard extremum seeking setup. The first one is that the problem is {\it time-varying}. In most of the classic extremum seeking literatures, the basic assumption is that the static input-output mapping is time invariant, i.e., \cite{Krstic2000,Tan2010,PID_ES} and the references therein. The time variation of the RPO problem may arise from the operation switching. The second one is {\it finite time horizon}. Unlike standard extremum seeking problem allowed to be solved in infinite time horizon, the RPO problem has to be solved with a finite duration due to its finite cycle duration. The third one, the most important one, is that the RPO problem has considerable {\it repetitiveness} thanks to the repeated operation mode. That makes room for improving transient performance by utilizing the repetitiveness, which {{motivates}} the introduction of iterative learning control (ILC). Finally, from the perspective of ILC, its standard assumption is that the tracking reference should be readily known \cite{Bristow}.}

In this paper, we propose a novel {{iterative learning control scheme}} based on {{extremum seeking}} to find the optimizer (optimal trajectory) of a time-varying mapping. Particularly, the contributions of the paper are threefold and summarized below.

First, {the proposed approach can be deemed as a novel extremum seeking scheme particular for time-varying optimizer tracking problem, because an analog memory has been introduced into the extremum seeking loop. That differs from the existing methods such as \cite{Wang,Guay_flat,Haring,Krstic2000sys,Krstic_book,Sahneh}.} There are two categories of similar approaches concerning the time-varying mappings reported in literatures. Wang and Krsti{\'c} introduced a detector to minimize the amplitude of stable limit cycle by tuning a controller parameter to a constant optimizer \cite{Wang}. Guay and his colleagues employed system flatness to parameterize all the variables by sine and cosine series; extremum seeking was used to steer the coefficients of the series to the optimizer \cite{Guay_flat}. Haring and his coworkers have developed a mean-over-perturbation-period filter to produce an estimate of the gradient for extremum seeking loop \cite{Haring}. The underlying assumption of all these methods is that the corresponding cost functions, although time varying, admit a constant optimizer, which is different from our discussion in this paper. As for the second type, 
Krsti{\'c} introduced a compensator for time-varying mappings which is structured by Wiener-Hammerstein models. This method requires the knowledge about the two time-varying blocks, which may restrict its applicability \cite{Krstic2000sys,Krstic_book}. An adaptive delay-based estimator was introduced to feed gradient estimates to extremum seeking loop by Sahneh and his coworkers. The extremum seeking loop is a cascade feedback in nature \cite{Sahneh}. Basically, this type of methods employed a fast decay ratio to suppress the tracking error. However, in some circumstances, these methods may not be able to yield satisfactory transient tracking performance. Fortunately, many time-varying systems exhibit certain repetitive behaviors \cite{Survey}, such as semiconductor manufacturing, pharmaceutical producing and the injection molding (c.f. Section \ref{PMFV}). Exploiting such repetitiveness provides potential for circumventing the aforementioned drawbacks. Actually, the works in \cite{Guay_flat, Haring} have already utilized such a feature and {{used it to formulate a periodic cost function.}} Iterative learning control, first proposed by Arimoto et. al. \cite{Arimoto}, is good at exploiting repetitiveness to improve tracking performance from iteration to iteration \cite{Bristow,Xu_survey,Survey,RDZ}. 

\looseness-1 Second, {the proposed approach can be categorized as an extension to ILC controller family, since it does not rely on the ``direction" of the ``feedbdack'' error information.} As mentioned before, the fundamental assumption adopted in most ILC literature is that the tracking reference must have been already available as priori knowledge \cite{Moore}, which renders the controller knowing the direction to steer. Within this paper, we do not need such an assumption; only the distance to the tracking reference should be known, i.e., the absolute value of tracking error. {For a simple example, let us consider the scalar system $y(t)=u(t)$ to track $u^*(t)$ with a proportional ILC control law $u_k(t)=u_{k-1}(t)+K[u_{k-1}(t)-u^*(t)]$. The corresponding error system $\tilde {u}_k(t)=(I+K)\tilde {u}_{k-1}(t)$ is convergent if and only if $|1+K|<1$. This can be achieved only when both $\text{sign}[u_{k-1}(t)-u^*(t)]$ and $|u_{k-1}(t)-u^*(t)|$ are known. If only $|u(t)-u^*(t)|$ is available, then $|1+K|<1$ cannot be ensured; thus, almost all the ILC approaches including norm optimal ILC \cite{Bristow}  will fail for this situation.} It is natural to come up with combining extremum seeking with ILC; let them {{collaborate with each other}}: ILC provides the past learning experience to extremum seeking to improve transient tracking performance; the ``direction'' information needed by ILC is {{given}} by extremum seeking. Utilizing extremum seeking to detect the ``direction" has been reported working quite well \cite{Unknown_direction}.

Third, {it can be shown that the proposed approach in nature turns out to be a new online infinite-dimensional integral-type (I-type) ILC control law, which has been analyzed with a new tool -- $k$-dependent contraction mapping.} Similar results are \cite{PILC,WangD,Saab,Saab2,Ouyang}. \cite{PILC} has studied the proportional-type (P-type) online ILC and derived an index bound regarding the ultimate tracking performance. Wang considered the sampling effect and input saturation issues in the offline P-type ILC, and implemented it experimentally \cite{WangD}. In \cite{Saab}, Saab investigated the offline P-type and D-type (derivative-type) ILC for the stochastic scenario, where a dynamic learning gain was adopted. Ref. \cite{Saab2} discussed the forgetting factor selection for a general offline ILC algorithm. Ouyang and his colleagues developed an online PD-type ILC for a class of input-affine nonlinear system, and also presented an ultimate bound of tracking performance \cite{Ouyang}. According to authors' knowledge, there is few papers contributing to online I-type ILC. Furthermore, we present more than an ultimate bound of tracking performance; the limit solution and its uniqueness are studied as well. More interestingly, an approximating system named \textit{modified Lie bracket system} has revealed that the proposed approach is essentially an online integral-type ILC with the approximation error as ``disturbance''. Based on that, it naturally extends the results on ILC to analyze the iterative learning extremum seeking system. We show that the system is uniformly bounded in terms of $k$  and converges to a set as $k$ goes to infinity. The particular set is a $\lambda$-norm ball, whose center is the limit solution of the associated ILC system, and radius can be controlled by the frequency of the sinusoid signal.

The rest of the paper is structured as follows: Section II gives the technical preliminaries about $\lambda$-norm and Lie bracket system; Section III formulates the problem; Section IV gives the analysis tool; Section V presents the main results; Section VI provides illustrative examples for the theory; a conclusion is drawn and an outlook is given in Section VII.

\textit{Notations:} $\mathbb{N}_{++}$ and $\mathbb{N}_0$ denote the set of positive integers excluding and including zero respectively. $\mathbb{M}_n$ is for all the matrices with dimensions $n\times n$. $C^n$ with $n\in\mathbb{N}_0$ stands for the set of $n$ times continuously differentiable functions and $C^\infty$ for the set of smooth functions. The gradient of a continuous function $f\in C^1:\mathbb{R}^n\to\mathbb{R}$ is $\nabla_x f(x)\triangleq \left[\frac{\partial f(x)}{\partial x_1},\dots,\frac{\partial f(x)}{\partial x_n}\right]^T$. Two vector fields $f,g:\mathbb{R}^n\times\mathbb{R}\to\mathbb{R}^n$ are twice continuously differentiable; their Lie bracket is defined as $[f,g](x,t)\triangleq \frac{\partial g(x,t)}{\partial x}f(x,t)-\frac{\partial f(x,t)}{\partial x}g(x,t)$. For a point $x\in \mathcal{X}$ and a set $\mathcal{S}\subset\mathcal{X}, x\notin\mathcal{S}$, the distance from $x$ to $\mathcal{S}$ is defined as $\text{dist}(x,\mathcal{S})=\inf_{s\in\mathcal{S}}\|x-s\|$. For two compact sets $\mathcal{X},\mathcal{Y}$ and $\mathcal{X}\subset\mathcal{Y}$, $\partial {\mathcal{Y}}$ for the boundary of $\mathcal{Y}$, the distance is defined as $\text{dist}(\mathcal{X},\mathcal{Y})=\min_{x\in\mathcal{X},y\in\partial\mathcal{Y}}\|x-y\|$. $\text{int }\mathcal{X}$ means the interior of set $\mathcal{X}$. %{{\underline{\color{blue}{The two defs are used in the proof of Lemma 3.2. }}}}

\section{Preliminaries}
\subsection{$\lambda$-norm}The $\lambda$-norm, introduced by Arimoto et. al. in 1984 \cite{Arimoto}, is a topological measure widely used to analyze the convergence of ILC control law \cite{lambda}. The formal definition of $\lambda$-norm is as follows.
\begin{definition}\cite{lambda}
The $\lambda$-norm of a function $f:[0,L]\to\mathbb{R}^n$ is
\begin{equation*}
\|f(\bullet)\|_\lambda=\max_{t\in[0,L]}\text{e}^{-\lambda t}\|f(t)\|_\infty
\end{equation*}
where $\|f(t)\|_\infty=\max_{1\le i\le n}|f_i(t)|$.
\end{definition}

From the definition, it is easy to see that
\begin{equation*}
\|f(\bullet)\|_\lambda\le\|f(\bullet)\|_C\le \text{e}^{\lambda L}\|f(\bullet)\|_\lambda
\end{equation*}
for positive $\lambda$, where $\|f(\bullet)\|_C=\max_{t\in[0,L]}\|f(t)\|_\infty$. This shows that the $\lambda$-norm is equivalent to the C-norm, which means the convergence with respect to $\lambda$-norm is still valid with respect to C-norm. Its advantage is that a non-monotonically converging sequence on C-norm can be monotonically converging on $\lambda$-norm for a properly chosen $\lambda$. 

\subsection{Lie bracket system} In the classic extremum seeking literatures, for example \cite{Krstic2000},  the behavior of the original extremum seeking system is analyzed by averaging. However, within this paper, an {{emerging analysis tool based on}} the Lie bracket approximation is going to be used to study the extremum seeking system. For an input-affine extremum seeking system 
\begin{equation*}
\dot{x}=b_1(x)\sqrt{\omega}u_1(\omega t)+b_2(x)\sqrt{\omega}u_2(\omega t)
\end{equation*}
with $\omega\in(0,\infty)$, its Lie bracket system is
\begin{equation*}
\dot{z}=\frac{1}{2}[b_1,b_2](z)
\end{equation*}
For instance, in a traditional ES system, $b_1(x)=1,b_2(x)=-\alpha f(x)$, $\alpha>0, f(x)\in C^2:\mathbb{R}^n\to\mathbb{R}$ admitting a local minimum, its Lie bracket system is $\dot{z}=-\alpha  \nabla f(z)/2$, which  clearly minimizes the cost function. {{Compared to standard averaging techniques, Lie bracket approximation techniques are easier to apply and provide appealing stability and convergence rate results. For details, please refer to \cite{Durr}.}}
%{{One advantage of using Lie bracket approximation}} is that the approximating system and the original system share the same rate of convergence, i.e., the {{exponential}} stability of the Lie bracket system implies the practical {{exponential}} stability of the original one. For details, please refer to \cite{Durr}.

\section{Problem statement}
{In this paper, we follow the standard setup for extremum seeking for time-varying mapping optimization problem in Ariyur and Krsti{\'c} (pp. 21) \cite{Krstic_book}. This parameterization can approximate a large class of vector function $F(x(t),t)$ admitting a quadratic time-varying minimizer $x^*(t)$. The problem can be formulated as follows: 
}{{at each time $t\in[0,L]$, we {{attempt to iteratively}} solve
\begin{equation}\label{eqn:opt} 
\min_{x(t)}F(x(t),t)
\end{equation}
where $x(t)\in\mathbb{R}^n$ and $F: \mathbb{R}^n\times[0,L]\rightarrow\mathbb{R}$,
\begin{equation}\label{eqn:map}
F(x(t),t)=f^*(t)+[x(t)-x^*(t)]^TQ[x(t)-x^*(t)]
\end{equation}
with $Q\in\mathbb{M}_n$ positive definite and $L>0$ is the finite time duration or period. For each time $t$, $x(t)=x^*(t)$ is the solution to the optimization problem (\ref{eqn:opt}) with the corresponding optimal value $f^*(t)$. For all time $t\in[0,L]$, all the $x^*(t)$ form a time-varying optimizer/minimizer trajectory (or called optimizer/minimizer) $x^*:[0,L]\to\mathbb{R}^n$; if $x^*(t)\equiv c$ for any $t$, $c$ a constant vector, then $x^*$ is named a constant optimizer. {{Similarly, the optimal value trajectory $f^*:[0,L]\to \mathbb{R}$ collects all the optimal value $f^*(t)$ for any $t$ within the interval $[0,L]$.}} $x(t)$ in (\ref{eqn:opt}) or its ensemble $x$ is named optimization variable or simply input. Obviously, the study on minimization problem {does not lose} any generality; it can be easily extended to the maximization problem by only altering the sign of the extremum seeking gain. \textit{The objective of extremum seeking is to iteratively steer {{the trajectory}} $x$ to the unknown optimizer trajectory $x^*$ over $[0,L]$. }}}
\begin{figure}[htbp]
\begin{center}
\includegraphics[width=8.2cm]{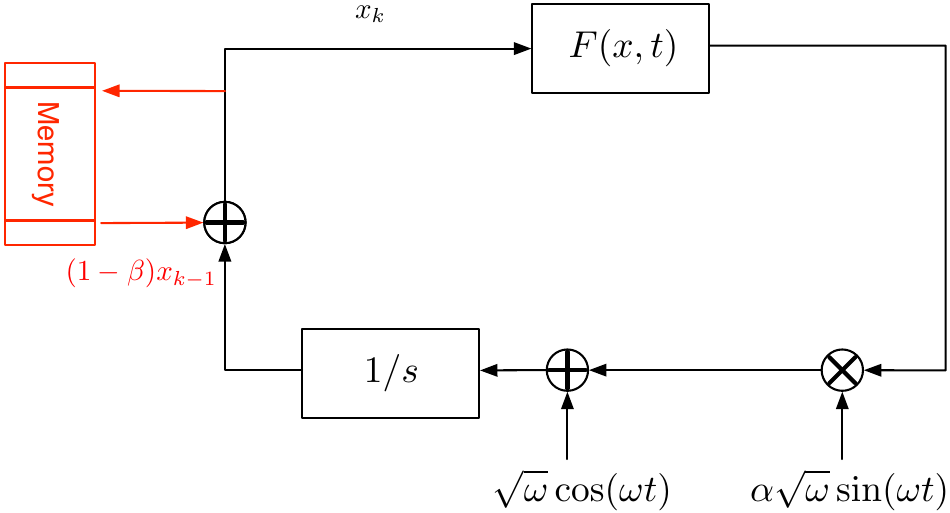}
\caption{Block diagram of iterative learning extremum seeking}
\label{fig1}
\end{center}
\end{figure}

\subsection{Classical approach}
The input $x$ in the conventional ES {{literature}}, is driven to {{a}} neighborhood of the optimizer trajectory $x^*$ by solving the following ordinary differential equation (ODE).
\begin{equation}\label{eqn:trad_es}
\dot{x}(t)=-\alpha F(x(t),t)\sqrt{\omega}\sin(\omega t)+\sqrt{\omega}\cos(\omega t)
\end{equation}
Here $\alpha\in(0,\infty)$ is a constant gain, {{and}} {{$\omega$ is the frequency of the perturbation signal.}} The solid-line part in Fig. \ref{fig1} presents the closed loop of this classic approach. For a fast changing optimizer trajectory, this approach may fail to achieve a satisfactory tracking within a finite duration.

\subsection{Main idea}
{{Because of}} the repeated operation mode of the repetitive process, there is no need to solve the tracking problem in only single iteration; instead, it can be solved in many iterations {{even infinite iterations}}. {{Since the previous input is a good approximation of the optimizer trajectory}}, it is quite handy and natural to modify the previous input to generate a new control input. Borrowing the ideas from ILC, we propose to solve the ES tracking problem by solving the following
ODE.
\begin{equation}\label{eqn:iles1}
\dot{x}_k(t)=(1-\beta)\dot{x}_{k-1}(t)-\alpha F(x_k(t),t)\sqrt{\omega}\sin(\omega t)+\sqrt{\omega}\cos(\omega t)
\end{equation}
The subscript $k\in\mathbb{N}_{++}$ indicates the iteration index; {{$\beta\in(0,1]$}} is the forgetting factor, which has been adopted in {{many ILC literature, for example \cite{Bristow,PILC}}}. {{Furthermore, $\dot{x}_0(t)$ and $x_0(t)$ are set to be zero.}} To keep the notation simple, $x_k$ will be used in the rest of paper instead of $x_k(t)$ without any ambiguity and the same for other similar variables, unless stated otherwise. 

Physically, (\ref{eqn:iles1}) introduces a memory storing the input $x_{k-1}$ of the previous batch into the extremum seeking loop as the red part shown in Fig. \ref{fig1}. The term $x_{k-1}$ is the feed-forward component, while the second and third terms in the right hand side of (\ref{eqn:iles1}) are the feedback component. It is intuitively understandable that the proposed method could result in a better performance than the conventional one (\ref{eqn:trad_es}), since the feed-forward term somehow can {{facilitate the tracking}}. Meanwhile, the mechanism is always feeding a new input into the system by insistently using feedback information to ``polish" $x_{k-1}$, a rough ``guess'' of the optimizer trajectory. Thus, it can be expected that the tracking performance would improve gradually as the rough ``guess" is becoming finer.% {{\underline{\color{blue}{Brief reasoning for why better.}}}}
\begin{remark}
It is noted that when $k=1$, (\ref{eqn:iles1}) becomes $\dot{x}_1=(1-\beta)\dot{x}_0-\alpha F(x_1,t)\sqrt{\omega}\sin(\omega t)+\sqrt{\omega}\cos(\omega t)$, which is exactly the standard extremum seeking, under the condition that $\dot{x}_0(t)$ and $x_0(t)$ are zero, or equivalently the memory is reset to zero.
\end{remark}

Since {{we are only interested in the system over the finite interval}}, there is no need to discuss its asymptotic stability along the time direction. The problems we are more interested in are under what condition $x_k$ will approach to a small neighborhood of $x^*$ when $k$ tends to infinity, i.e., {{$\|x_k-x^*\|_C<D, k\to\infty$}} and what determines $D$. 
 
\section{Preliminary results}
{Within the section, an analysis tool named {\it $k$-dependent contraction mapping} is developed in Theorem \ref{thm-4-1}. It is general, because it can be applied {{to an infinite dimensional setting}}, and it also lays the {{foundation}} for analyzing ILES.}

 An operator $T$ between two real linear spaces $\mathcal{X}$ and $\mathcal{Y}$ is called a {\it linear mapping} or {\it linear operator} if $T(\lambda x+\mu y)=\lambda T(x)+\mu T(y)$ for all $\lambda,\mu\in\mathbb{R}$ and $x,y\in\mathcal{X}$. There is a norm $\|\bullet\|$ defined on $\mathcal{X}$ and $\mathcal{Y}$. Then, a linear mapping is {\it bounded} if there exists a constant $M\ge 0$ such that
\begin{equation*}
\|T(x)\|\le M\|x\| \text{ for all  } x\in\mathcal{X}.
\end{equation*}
The sets consisted of all these bounded linear mapping $T$ is denoted by $\mathsf{B}(\mathcal{X},\mathcal{Y})$. We write $\mathsf{B}(\mathcal{X})$ for short if the domain and range spaces are the same. Moreover, the {{operator norm}} is defined as
\begin{equation*}
\|T\|=\sup_{x\neq 0}\frac{\|T(x)\|}{\|x\|}.
\end{equation*}
\begin{definition}[Uniformly convergence (pp.109, \cite{Hunter})] 
If $\{T_k\}$ is a sequence of mappings in $\mathsf{B}(\mathcal{X},\mathcal{Y})$ and 
\begin{equation*}
\lim_{k\to\infty}\|T_k-T\|=0
\end{equation*}
for some $T\in\mathsf{B}(\mathcal{X},\mathcal{Y})$, then we say that $T_k$ converges uniformly to $T$.
\end{definition}

Note that $\|T_k(x)-T(x)\|\le\|T_k-T\|\|x\|$. Given that $\|x\|$ is bounded, $\lim_{k\to\infty}\|T_k-T\|=0$ implies that $\lim_{k\to\infty}T_k(x)=T(x)$. In other words, uniform convergence implies {\it strong convergence}. 

\begin{theorem}[$k$-dependent contraction mapping]\label{thm-4-1}
Let $S$ be a closed and bounded subset of a Banach space $\mathcal{X}$. If a mapping sequence $\{T_k\}$ satisfies
\begin{itemize}
\item
C1) for every $k$, $T_k\in\mathsf{B}(S)$;
\item
C2) $\|T_k(x)-T_k(y)\|\le\rho\|x-y\|,\forall x,y\in S$, for a universal \emph{contraction mapping ratio} $\rho\in[0,1)$, denoting $\mathcal{T}$ as the set {consisting} of all such $T_k$;
\item
C3) $T_k$ converges uniformly to $T_\infty$ as $k\to\infty$, for some $T_\infty\in\mathcal{T}$,
\end{itemize}
then $T_k$ is called a $k$-{\it dependent contraction mapping} on $S$. Furthermore, 
\begin{itemize}
\item
there exists a unique solution $x^\star\in S$ satisfying $x^\star=T_\infty(x^\star)$;
\item
$x^\star$ is independent of the initial value $x_1\in S$; {$x_1$ is arbitrarily selected in $S$}.
\end{itemize} 
\end{theorem}
The proof is provided in Appendix.
 
\section{Main results}
{In Section \ref{V_MLB}, Lemma \ref{lem1} shows the existence of contraction mapping ratio as to {\it modified Lie bracket (MLB) system}; it verifies the satisfaction of C2 in Theorem \ref{thm-4-1}. It also shows that the MLB system of the iterative learning extremum seeking is an I-type ILC. In Section \ref{IVA}, the approximation-error-free system (I-type ILC) is studied. Lemma \ref{lem_ilc1} shows the existence of an invariant set for the mapping; Theorem \ref{thm_ilc1} studies the approximation-error-free system with $\beta\neq 0$; Theorem \ref{thm_ilc2} discusses the case $\beta=0$. In Section \ref{IVB}, Lemma \ref{lem2} shows the original system (\ref{eqn:iles1}) and approximate system (\ref{eqn:mlbs}) can be close enough; Theorem \ref{thm} proves the uniform boundedness of the original system (\ref{eqn:mlbs}); Theorem \ref{thm_4_2} studies the property of the limit solution.
}

\subsection{MLB system}\label{V_MLB}
Note that if (\ref{eqn:trad_es}) iterates itself to different {{$k$}}, the corresponding coefficients on $\cos(\omega t)$ will be different. On the other hand, the Lie bracket approximation is only valid for systems {{with}} respect to $t$ not to both $k$ and $t$. Thus, we cannot directly employ it but have to tailor it accordingly. We propose to use the MLB system to approximate the original one.

%{{\underline{\color{blue}{Reason for introducing $\gamma_k$.}}}} 
\begin{equation}\label{eqn:mlbs}
\begin{aligned}
\dot{z}_k=&(1-\beta)\dot{x}_{k-1}+\frac{\gamma_k}{2}[1,-\alpha F](z_k,t)\\
=&(1-\beta)\dot{x}_{k-1}-\frac{\alpha\gamma_k}{2}\nabla_z F(z_k,t)
\end{aligned}
\end{equation}
Here the only modification is the introduction of $\gamma_k$. $\gamma_k$ is a compensating parameter and only related to the forgetting factor $\beta$ as defined below.
\begin{equation}\label{eqn:gamma_k}
\gamma_k=\frac{1-(1-\beta)^k}{\beta}
\end{equation}
It is evident that $\{\gamma_k\}$ is a monotonically increasing sequence and $1\le\gamma_k\le1/\beta$. The first equality holds if and only if $k=1$, which implies that the MLB system (\ref{eqn:mlbs}) reduces to the traditional Lie bracket system when $k=1$. {The reason to introduce $\gamma_k$ is to compensate the gradient mismatch between $x_k$ and $z_k$, which arises from the accumulating effects of successive iterations. To see this, please refer to (\ref{eqn:iter}) and (\ref{eqn:xz1}).} The term $\dot{x}_{k-1}$ is a feed-forward term and only a function of time with respect to the current iteration $k$, {{and does not contain any sinusoid term}}. Thus, it should not be included into the Lie bracket according to the theory in \cite{Durr}.
\begin{remark}
Observing (\ref{eqn:mlbs}), the MLB system (\ref{eqn:mlbs}) is unlike the conventional Lie bracket system, because of the existence of the {time derivative} of $x_{k-1}$ rather than being an independent system of itself. Inserting (\ref{eqn:map}) into (\ref{eqn:mlbs}) and denoting
\begin{equation*}
\Gamma_k=\frac{\alpha\gamma_kQ}{2}=\frac{[1-(1-\beta)^k]\alpha Q}{2\beta}
\end{equation*}
\end{remark}
we can rewrite (\ref{eqn:mlbs}) as
\begin{equation}\label{eqn:mlbs1}
\dot{z}_k=(1-\beta)\dot{z}_{k-1}-\Gamma_k(z_k-x^*)+(1-\beta)(\dot{x}_{k-1}-\dot{z}_{k-1}).
\end{equation}
{It should be noted that $x^*$ in {{(\ref{eqn:mlbs1})}} is only used for conceptual analysis; in fact, we do not require the knowledge of $x^*$ in implementation, because only (\ref{eqn:iles1}) is implemented in practice.}

Before presenting our main results, we impose the following assumptions.
\begin{itemize}
\item
A1) The time-varying optimizer trajectory $x^*\in C^1: [0,L]\to\mathbb{R}^n$ and optimal value trajectory $f^*\in C^1: [0,L]\to\mathbb{R}$;
\item
A2) The initial condition of each iteration of (\ref{eqn:iles1}) is identical and equal to zero, i.e. $x_k(0)=x(0)=0, \forall k$; so is the MLB system (\ref{eqn:mlbs}), i.e. $z_k(0)=x(0)=0$ for all $k$; 
\item
A3) Assume that {{$\|Q\|$}} is bounded and $Q\ge\delta I$, where $\delta$ is a known positive real number.
\end{itemize}

\begin{remark}
A1 is assumed to ensure the existence of $F(x,t)$'s first-order partial derivative $\frac{\partial F(x,t)}{\partial t}$, which is required by integration by parts. A2 is a {{common}} assumption in majority of ILC literature to simplify derivations \cite{Bristow}, called \textit{identical initial condition (i.i.c)}. The second part of A2 is assumed to let MLB system be in accordance with the original system. A3 means that we do not require exact knowledge about $Q$, the Hessian matrix, but requires a lower bound of $Q$ for the minimum case, since a known $Q$ means having a precise knowledge about the plant, which is generally impossible {{in practice}}. $\delta$ is only required by the following conceptual analysis, does not restrict our method's applicability.
\end{remark}

Taking integration on both sides of (\ref{eqn:mlbs1}), it turns to be
\begin{equation}\label{eqn:i-ilc}
z_k=(1-\beta)z_{k-1}-\Gamma_k\int_0^t(z_k-x^*)ds+(1-\beta)(x_{k-1}-z_{k-1})
\end{equation}
The $z_k$ can be interpreted as control input with $z_k-x^*$ as tracking error. Then, the rewriting above clearly shows {{that it is in the form}} of integral-type ILC online (feedback) control {{\cite{Bristow}}} with the approximation error, i.e., $(1-\beta)(x_{k-1}-z_{k-1})$.%\hfill{\scriptsize$\blacksquare$} 

Note that {{(\ref{eqn:mlbs})}} is a linear ordinary differential equation, and we can write down its explicit solution.
\begin{equation*}
z_k=\text{e}^{-\Gamma_k t}z(0)+\int_0^t\text{e}^{-\Gamma_k(t-s)}[\Gamma_kx^*+(1-\beta)\dot{x}_{k-1}]ds
\end{equation*}
It follows from $z(0)=0$ and integrating by parts that
\begin{equation*}
\begin{aligned}
z_k=&(1-\beta)\text{e}^{-\Gamma_kt}\left[\left.\text{e}^{\Gamma_ks}x_{k-1}\right|_0^t-\int_0^t\text{e}^{\Gamma_ks}\Gamma_kx_{k-1}ds\right]\\
&+\text{e}^{-\Gamma_kt}\int_0^t\text{e}^{-\Gamma_ks}\Gamma_kx^*ds
\end{aligned}
\end{equation*}
Because $x_{k-1}(0)=0$, we have
\begin{equation}\label{eqn:exp_sol}
z_k=(1-\beta)x_{k-1}+\text{e}^{-\Gamma_kt}\int_0^t\text{e}^{\Gamma_ks}\Gamma_k[x^*-(1-\beta)x_{k-1}]ds
\end{equation}
Now we define the tracking error of the MLB system (\ref{eqn:mlbs}) as $y_k\triangleq z_k-x^*$; the error system is
\begin{equation}\label{eqn:yk1}
y_k=T_k\left(y_{k-1}+x_{k-1}-z_{k-1}-\frac{\beta}{1-\beta}x^*\right)
\end{equation}
where $T_k$ is the mapping as follows.%:[0,L]\times[0,L]\to\mathbb{R}^n$
\begin{equation}\label{eqn:tk}
T_k(x)(t)=(1-\beta)x(t)-(1-\beta)\text{e}^{-\Gamma_kt}\int_0^t\text{e}^{\Gamma_ks}\Gamma_kx(s)ds
\end{equation}
For the sake of simple notation, we will write (\ref{eqn:tk}) for short as $T_k(x)=(1-\beta)x-(1-\beta)\text{e}^{-\Gamma_kt}\int_0^t\text{e}^{\Gamma_ks}\Gamma_kxds$. Note that the mapping above is $k$-dependent because $\Gamma_k$ is $k$-dependent. Now we will give the result of contraction mapping for $k$-dependent case, which differs from pp. 655 \cite{Khalil} (where only $k$-invariant mappings are studied). 

Define a Banach space ${\cal{X}}=C[0,L]$ and a closed and bounded set %(Example 13.18 in \cite{Hunter1})
\begin{equation}\label{eqn:set_S}
S=\{y\in{\cal{X}}|\|y\|_\lambda\le D_2\}
\end{equation}
{Let $y$ be the tracking error of MLB system defined as $y=z-x^*$.} Then (\ref{eqn:set_S}) implies that $z_k$ is contained in $\mathcal{S}_z=\{z\in\mathcal{X}|\|z-x^*\|_\lambda\le D_2\}$. The lemma below shows the conditions to fulfill C2 for $T_k$.

\begin{lemma}\label{lem1}
Consider the mapping (\ref{eqn:tk}) and let A1-A3 hold. For arbitrary $\beta\in(0,1)$, there is a $\lambda_0$ such that $\|T_k(x)-T_k(y)\|_\lambda\le\rho\|x-y\|_\lambda$ for any $x,y\in\mathcal{X}$, {{$\rho\in(0,1)$}} for every $\lambda\in(\lambda_0,\infty)$. Moreover, $\rho$ and $\lambda_0$ are independent of $k$.
\end{lemma}
%Or equivalently, there exists a $\rho\in(0,1)$ such that $\|T_k(x)-T_k(y)\|_\lambda\le\rho\|x-y\|_\lambda$.}
%Conversely, given $\lambda$ there also exists a $\beta_0$ such that $T_k$ is a contraction mapping in terms of $\lambda$-norm for every $\lambda\in(\beta_0,1)$.
%\begin{IEEEproof}
%For details of the proof, please refer to Appendix \ref{pf1}.
%\end{IEEEproof}

%{{\color{blue}{A3 should hold, because it is used to describe how large $\lambda$ should be. And I have make it clear that A3 is only needed by the conceptual proof not by practice.}}}

\begin{remark}\label{rem1}
Many ILC literatures do not use such a forgetting factor, since its existence will compromise the zero tracking error (perfect tracking) iteration-wisely \cite{PILC}. However, the above lemma suggests otherwise in our case. It is necessary to introduce such a forgetting factor $\beta$ to ensure the existence of the $k$-dependent contraction mapping, since that the perfect tracking cannot be achieved due to the existence of dither signal. If without $\beta$, this undesired effect would keep accumulating and the overall performance would deteriorate rapidly.
\end{remark}
\begin{remark}
If the C-norm is used instead of $\lambda$-norm, or equivalently $\lambda=0$, then according to the proof of Lemma \ref{lem1}, $\beta$ has to be greater than $2-\sqrt{2}$. It suggests that $\lambda$-norm somehow enlarges the feasible basin of $\beta$.
\end{remark}
\begin{remark}
From (\ref{eqn:tk}), it is evident that the mapping $T_k$ is a linear mapping. Thus, we can rewrite (\ref{eqn:yk1}) as
\begin{equation}\label{eqn:yk2}
y_k=T_k(y_{k-1})+T_k(x_{k-1}-z_{k-1})+T_k\left(-\frac{\beta}{1-\beta}x^*\right)
\end{equation}
In (\ref{eqn:yk2}), it can be seen that the second and third terms on the right hand side play a role like ``disturbances''; the second one is caused by the approximation, while the third is caused by the forgetting factor. The two ``disturbances'' are still different: $T_k(x_{k-1}-z_{k-1})$ is a persistently active noise, while $T_k[-\beta x^*/(1-\beta)]$ is just a constant offset. 
\end{remark}
%Recalling the fix-point theorem (cf. pp. 655 \cite{Khalil}), there are no such ``disturbance'' terms; thus, we can always restrict the input and output of our mapping in a closed set, i.e., $y_k,y_{k-1}\in S$ (ignoring the disturbance) in our particular case. From this point, it seems that we have to show $y_{k-1}+x_{k-1}-z_{k-1}-\beta x^*/(1-\beta)$ within $S$ on the way towards the conclusion that $T_k$ is a contraction mapping. In fact, we do not need to do so, since Lemma \ref{lem1} holds for $\mathcal{X}$ the whole. $T_k$ applying on $x_{k-1}-z_{k-1}$ and $\beta x^*/(1-\beta)$ is only for short hand. Moreover, owing to ``disturbance'', we will not show $\tilde{y_k}=T_k(\tilde y_{k-1}), y_1=\tilde y_1$ within $S$ towards claiming that $T_k$ is a contraction mapping; but instead, we will prove the stronger version of that, namely, $y_k$ in (\ref{eqn:yk1}) within $S$ for the worst case that ``disturbance'' is always driving the trajectory away. 

\subsection{{Approximation-error-free system\label{IVA}}}
{{Within this subsection, we temporarily assume that the approximation error, i.e. $x_k-z_k$ is zero. Thus, we will only study the behavior of the approximation-error-free MLB system (\ref{eqn:i-ilc}), i.e.,}}
%We are ready to present our contribution on ILC, which studies the behavior of the disturbance-free MLB system (\ref{eqn:i-ilc}), i.e., 
\begin{equation}\label{eqn:i-type-ilc}
z_k=(1-\beta)z_{k-1}-\Gamma_k\int_0^t(z_k-x^*)ds.
\end{equation}
{{Note that (\ref{eqn:i-type-ilc}) can be interpreted as a static system $z=u$ regulated by ILC control law $u_k=(1-\beta)u_{k-1}-\Gamma_k\int_0^te_k ds$ with $e_k$ being the tracking error. Since (\ref{eqn:i-type-ilc}) represents the dominant dynamics of (\ref{eqn:i-ilc}), studying (\ref{eqn:i-type-ilc}) is also helpful to understand iterative learning extremum seeking. The rest of this subsection is divided into two parts according to different values of forgetting factor $\beta$. For $\beta\neq 0$, the $k$-dependent contraction mapping will be used to study the dynamics of (\ref{eqn:i-type-ilc}), while for $\beta=0$, a Lyapunov-like argument helps to understand (\ref{eqn:i-type-ilc}).}}

It is equivalent to study the following error system instead of (\ref{eqn:i-type-ilc}).%Equivalently, this corresponds to be the error system
\begin{equation}\label{eqn:yk3}
y_k=T_k(y_{k-1})+T_k\left(-\frac{\beta}{1-\beta}x^*\right)
\end{equation}
%It should be mentioned that in this section, we impose that $x^*$ is known. %For a given $x\in S$, $T_k(x)$ is a continuous function of $\gamma_k$ and $\gamma_k\to 1/\beta$ as $k\to\infty$. Thus, there exists a $T_\infty$ for $T_k$ sequence. Consequently, C3 is fulfilled. Lemma \ref{lem1} suggests that C2 is satisfied. It remains to show C1. The following lemma demonstrates a stronger version of C1: with the ``offset'', i.e., $T_k[-\beta x^*/(1-\beta)]$, $y_k$ is uniformly bounded, namely, $D_2$ exists.
 Define the mapping
\begin{equation}\label{eqn:gk}
G_k(x)\triangleq T_k(x)+T_k\left(-\frac{\beta}{1-\beta}x^*\right).
\end{equation}
The goal is to show $G_k$ is a $k$-dependent contraction mapping so that the convergence of $y_k$ can be concluded. It is easy to verify that $G_k$ satisfies C2, since that $G_k(x)-G_k(y)=T_k(x)-T_k(y)$ and $T_k(x)-T_k(y)$ fulfills C2. The following lemma gives a sufficient condition that $G_k$ fulfills C1.

\begin{lemma}\label{lem_ilc1}
Consider the mapping in (\ref{eqn:gk}) and let A1-A3 be satisfied. Given {{$\lambda>0, \beta\in(0,1), \rho\in(0,1)$}}, if $D_2$ in (\ref{eqn:set_S}) satisfies
\begin{equation*}
D_2\ge\max\{D_0,D^*\},
\end{equation*}
then $G_k$ maps $S$ into $S$. $D^*$ is defined as
\begin{equation*}
D^*=\frac{\beta\rho}{(1-\beta)(1-\rho)}\|x^*\|_\lambda
\end{equation*}
{{and}} $D_0=\|y_1\|_\lambda$, $y_1$ is defined by
\begin{equation*}
y_1=-x^*-\Gamma_1\int_0^ty_1ds
\end{equation*}
It is simply executing (\ref{eqn:yk3}) on $k=1$. 
\end{lemma}

%\begin{IEEEproof}
%For details of the proof, please refer to Appendix \ref{pf_ilc1}.
%\end{IEEEproof}

Lemma \ref{lem_ilc1} not only presents a sufficient condition for $G_k$ satisfying C1, but also states that $\|y_k\|_\lambda$ is uniformly bounded. Obviously, we can offer more than that, i.e., convergence and uniqueness of the limit solution, by invoking the $k$-dependent contraction mapping theorem.
\begin{theorem}\label{thm_ilc1}
Consider the mapping in (\ref{eqn:gk}), {$\beta\neq 0$} and let A1-A3 hold. If $D_2\ge\max\{D_0,D^*\}$ as in Lemma \ref{lem_ilc1}, then there exists a unique limit solution {{$y_k$ {{of}} (\ref{eqn:yk3})}} as $k$ tends to infinity. Moreover, $\lim_{k\to\infty}y_k=y_\infty$. $y_\infty$ is defined as the solution to the following equation
\begin{equation}\label{eqn:y_inf}
y_\infty=T_\infty\left(y_\infty-\frac{\beta}{1-\beta}x^*\right).
\end{equation}
$T_\infty$ is the limit of $T_k$ as $k$ tends to infinity, i.e.,
\begin{equation}\label{eqn:t_inf}
T_\infty(x)=(1-\beta)x-(1-\beta)\text{e}^{-\Gamma_\infty t}\int_0^t\text{e}^{\Gamma_\infty s}\Gamma_\infty xds
\end{equation}
{{and}} $\Gamma_\infty=\frac{\alpha Q}{2\beta}$ with {$\beta\neq 0$}.
\end{theorem}

\begin{remark}
Theorem \ref{thm_ilc1} shows that the trajectory of the approximation-error-free MLB system (\ref{eqn:i-type-ilc}) will ultimately converge to a fixed trajectory, which is parameterized by $\beta$ and $x^*$. From (\ref{eqn:t_inf}), it is clear that $0$ is a solution to the equation $x=T_\infty(x)$. From (\ref{eqn:y_inf}), one can see that $y_\infty$ will approach 0 if $\beta\to0$. However, will $y_\infty$ be $0$ if $\beta=0$? From Lemma \ref{lem1}, it is known that the contraction mapping method will fail when $\beta=0$. Thus, we use another way to prove the claim in the theorem below.%Taking partial derivative on $T_\infty$ with respect to $x$, we also have
\end{remark}

Before presenting the theorem, however, we have to make a remark on this particular case $\beta=0$. In the following, we will abandon $\Gamma_k$ but use a constant gain $\Gamma$ (positive definite matrix, not necessary restricting to $\alpha Q/2$) instead. The reasons for doing so are twofold. First, (\ref{eqn:gamma_k}) suggests that $\Gamma_k$ will be ill-defined, since {{$\gamma_k$ will be not defined when $\beta=0$}}. Second, as Remark \ref{rem1} states, the motivation of introducing $\beta$ is to handle the approximation error, and the gain $\Gamma_k$ becomes $k$-varying because of $\beta$; now we are hereby dealing with approximation-error-free control system.   

{{
\begin{theorem}\label{thm_ilc2}
Consider (\ref{eqn:i-type-ilc}) with $\beta=0$, that is $z_k$ is defined by the following formula
\begin{equation}\label{eqn:mlb0z}
\dot{z}_k=\dot{z}_{k-1}-\Gamma(z_k-x^*).
\end{equation}
If A1-A3 hold and $x^*(0)=0$, then
\begin{equation*}
z_k\to x^* \text{ almost everywhere as } k\to \infty.
\end{equation*}
\end{theorem}

\begin{remark}
Theorem \ref{thm_ilc2} gives a weaker result than Theorem \ref{thm_ilc1}, since it can converge to $0$ except on a set whose measure is $0$ and the uniqueness is lost. The result is quite understandable from a perspective of Laplace transform. Taking Laplace transform on both sides of (\ref{eqn:thm4-2-1}), we have
\begin{equation*}
Y_k(s)=\frac{s}{s+\Gamma}Y_{k-1}(s)
\end{equation*}
The modulus of the gain is less than 1, i.e., $|s/(s+\Gamma)|<1$. But it tends to $1$ as $s\to\infty$. It suggests that this algorithm have a weaker decaying effect on high-frequency signal. It coincides with the result.
\end{remark}
}}

\subsection{{Iterative learning extremum seeking control}}\label{IVB}
We follow {{a}} similar idea {{as outlined above}} to study the dynamics of ILES. ILES is an ILC control policy with the approximation error as ``disturbance''. Due to the ``disturbance'', the unique limit solution cannot be achieved. Therefore, ILES converging to a set will be shown instead, if the ``disturbance'' is bounded.

According to (\ref{eqn:yk2}), we define a new operator $H_k$ as
\begin{equation*}
H_k(x)\triangleq T_k(x)+T_k[-\beta x^*/(1-\beta)]+T_k(x_{k-1}-z_{k-1})
\end{equation*}
Since $T_k$ is a bounded operator, so is $H_k$ provided that $x_{k-1}-z_{k-1}$ is bounded. It is obvious that $H_k(x)-H_k(y)=T_k(x)-T_k(y)$; it is easy to verify that $H_k$ satisfies C2. Following the procedures we did in Section IV.A, we are going to show that how to properly design $D_2$ to ensure $H_k$ maps $S$ into $S$.

The following lemma lays the foundation of {{inductive}} arguments towards that conclusion. It basically means that {{for any $k$ we can always ensure that $x_k$ is in a invariant set given that its MLB system $z_k$ is in the interior of the invariant set, if $x_1,x_2,\dots,x_{k-1}$ are all in that set.}} This can be achieved by selecting a frequency $\omega$ larger than a threshold $\omega_0$, which is independent of $k$ and $t$.
\begin{lemma}\label{lem2}
Let A1-A3 be satisfied. Given an integer $k_0$, considering $x_k$ in (\ref{eqn:iles1}), suppose that {{$x_k$ for $k=1,\dots, k_0-1$ are}} uniformly in $t$ contained within a compact set $\mathcal{S}\subseteq\mathbb{R}^n$ and $0\in\text{int }\mathcal{S}$. If $z_{k_0}$ in (\ref{eqn:mlbs}) is contained in $\mathcal{K}\subset \text{int }\mathcal{S}$, $0\in\mathcal{K}$, then there exists a $\omega_0\in(0,+\infty)$ {{such that for every ${{\omega}}\in(\omega_0,\infty)$,}} $x_{k_0}(t)$ is uniformly in $t$ contained within $\mathcal{S}$ as well for any $\beta\in(0,1)$. Moreover, $\text{dist}(\mathcal{K},\mathcal{S})$ can be made arbitrarily small by selecting a sufficiently large $\omega$.%{\color{blue}{It is the result of A2. It lets us not worry about the initial condition of each iteration.}}
\end{lemma}
\begin{IEEEproof}
This proof uses the similar arguments as Theorem 1 in \cite{Durr}, but tailors them for the iterative case. For details of the proof, please refer to Appendix.
\end{IEEEproof}

The lemma above also indicates that the approximation error $\|x_k-z_k\|_\lambda$ is uniformly bounded; it can be arbitrarily small by selecting a sufficiently large $\omega$. 
%To this end, we are ready to answer the question quantitively how close the original system can track the time varying optimizer by applying the extremum seeking law (\ref{eqn:iles1}). The question consists of two sub-questions: first, how close the MLB system (\ref{eqn:mlbs}) can follow the desired optimizer; second, whether the MLB system (\ref{eqn:mlbs}) can uniformly approximate the original extremum seeking system in a certain neighborhood, i.e., $\|x_{k}-z_{k}\|_C<D_1$ for all $k$. Both sub-questions are answered in the following theorem.

\begin{theorem}\label{thm}
Consider (\ref{eqn:yk2}) and let A1-A3 hold. Given $\lambda$ and $\rho\in(0,1)$, there exists a $\omega_0\in(0,+\infty)$ such that for every $\omega\in(\omega_0,\infty)$, $H_k$ maps $S$ into $S$, if $D_2$ in (\ref{eqn:set_S}) satisfies
\begin{equation*}
D_2\ge \max\{D_0,D^\star\}
\end{equation*}
where $D_0$ is defined in Lemma \ref{lem_ilc1} and 
\begin{equation*}
D^\star=\frac{\rho}{1-\rho}\left(D_1+\frac{\beta}{1-\beta}\|x^*\|_\lambda\right).
\end{equation*}
$D_1$ is the uniform bound of $\|x_k-z_k\|_\lambda$ for an arbitrary $k$.
\end{theorem}

Theorem \ref{thm} suggests that the uniform bound of tracking error of the MLB system (\ref{eqn:mlbs}) -- $D_2$, can be made small by reducing the approximation error ($D_1$), i.e., employing a sufficiently large $\omega$. 

Since $x_k-z_k$ is consistently varying, it is impossible to show that $H_k$ satisfies C3. Therefore, we cannot conclude the unique limit solution; however, we can show that $y_k$ converges to a $\lambda$-norm ball. 
%It is interesting to know what will happen for the ideal case $\omega\to\infty$, since $x_k-z_k$ will disappear in (\ref{eqn:yk1}). In fact, we can conclude the existence and uniqueness of $\lim_{k\to\infty}y_k$ rather than $\lim_{k\to\infty}\|y_k\|_\lambda$ for $\omega\to\infty$. %{\color{blue}{As I mentioned before, we cannot conclude the existence of $\lim_{k\to\infty}x_k-z_k$ even for a sufficiently large $\omega$, which make us incapable to further claim something on $\lim_{k\to\infty}y_k$.}}
%Now, we are going to discuss the existence of $\lim_{k\to\infty}y_k$ for the extreme case $\omega\to\infty$.

%{\color{blue}{The following uses fix point theorem but not exactly the same.}}
\begin{theorem}\label{thm_4_2}
Consider (\ref{eqn:yk2}) and let A1-A3 be satisfied. Given $\lambda$, {{$\beta\in(0,1)$}} and $\rho\in(0,1)$, there exists a $\omega_0$ such that for every $\omega\in(\omega_0,\infty)$, {{$y_k$ in (\ref{eqn:yk2})}} will converge to a set $\mathcal{Y}$ as $k$ tends to infinity. Furthermore, 
\begin{equation}\label{eqn:set_y}
\mathcal{Y}=\{y\in S|\|y-y_\infty\|_\lambda\le D_y\}
\end{equation}
Here $y_\infty$ is defined in (\ref{eqn:y_inf}), $D_y=\rho D_1/(1-\rho)$; $D_1$ is the uniform bound of $\|x_k-z_k\|_\lambda$ for an arbitrary $k$. 
\end{theorem}

%\begin{IEEEproof}
%The proof is in Appendix \ref{pf_4_2}.
%\end{IEEEproof}

\begin{remark}
Theorem \ref{thm_4_2} suggests that we cannot achieve ``perfect tracking'' or a fixed limit trajectory unlike many ILC control laws, due to the existence of dither signals (sinusoid signals). However, according to Lemma \ref{lem2}, $D_1$ can be made arbitrarily small for a sufficiently large $\omega$. Thus, so is $D_y$, since $D_y$ is proportional to $D_1$; that means we can make the ultimate trajectory be as close to a fixed limit trajectory as one wishes by selecting a sufficiently large frequency. In the mean time, we can also let the fixed trajectory be as close to zero as possible by having a small enough forgetting factor $\beta$.
\end{remark}

\section{Illustrative {{examples}}}
\subsection{Approximation-error-free system}
{{{{In order to illustrate the results for the approximation-error-free system}} (Theorem \ref{thm_ilc1} and Theorem \ref{thm_ilc2}), we present the following numerical example.}} Consider the problem of tracking the following reference using (\ref{eqn:i-type-ilc}).
\begin{equation*}
x^*(t)=-\sin\left(\frac{\pi}{20}t\right)
\end{equation*}

\begin{figure}[htbp]
                \includegraphics[width=0.42\textwidth]{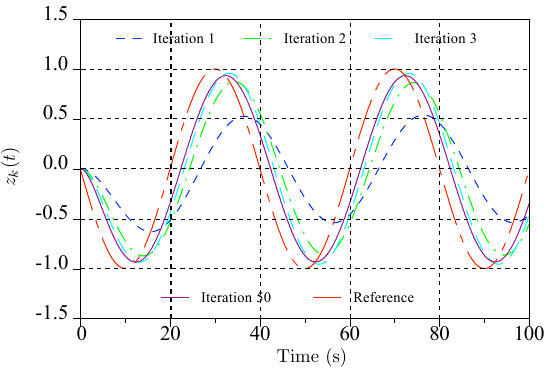}
                \caption{Approximation-error-free system with forgetting factor $\beta=0.5$}
                \label{ILC1}
\end{figure}

\begin{figure}[htbp]
                \includegraphics[width=0.42\textwidth]{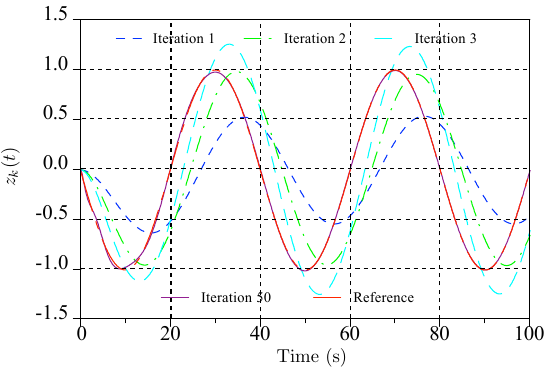}
                \caption{Approximation-error-free system without forgetting factor ($\beta=0$)}
                \label{ILC2}
\end{figure}

Fig. \ref{ILC1} shows the trajectory evolution versus iteration for the approximation-error-free system with a forgetting factor $\beta=0.5$. It clearly verifies that the trajectory ($z_k$) will converge to a fixed trajectory, but the trajectory has a gap with the tracking reference. Fig. \ref{ILC2} demonstrates the evolution for the approximation-error-free system without the forgetting factor $\beta$. It is shown that the trajectory will converge to the reference ultimately, thus ``perfect tracking'' achieved. Meanwhile, it {illustrates} the result of Theorem \ref{thm_ilc2}. {{It also suggests that there exists a tradeoff between convergence rate and tracking {{error}}.}}

\subsection{ILES}
\begin{figure*}[htbp]
\begin{subfigure}[b]{0.24\textwidth}
                \includegraphics[width=\textwidth]{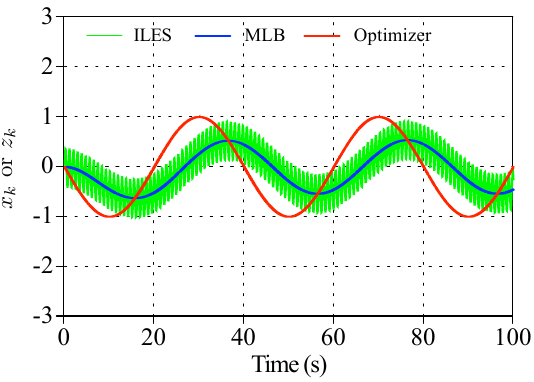}
                \subcaption{\scriptsize$1^\text{st}$ iteration with $\omega=7 {{\text{ rad/s}}}$}
                \label{fig2}
\end{subfigure}
\begin{subfigure}[b]{0.24\textwidth}
                \includegraphics[width=\textwidth]{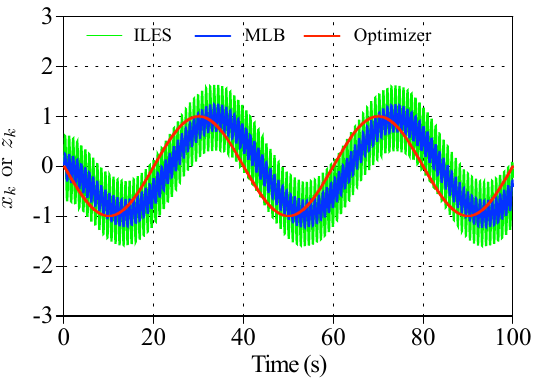}
                \subcaption{\scriptsize$2^\text{nd}$ iteration with $\omega=7 {{\text{ rad/s}}}$}
                \label{fig3}
\end{subfigure}
\begin{subfigure}[b]{0.24\textwidth}
                \includegraphics[width=\textwidth]{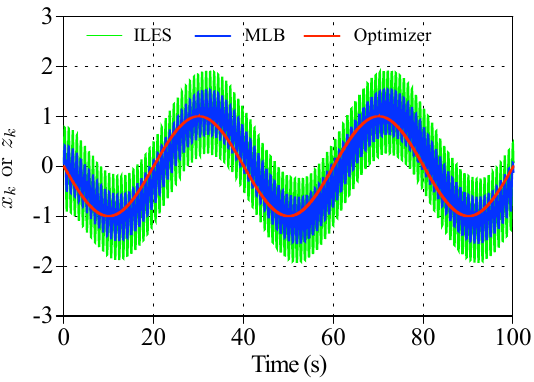}
                \subcaption{\scriptsize$3^\text{rd}$ iteration with $\omega=7 {{\text{ rad/s}}}$}
                \label{fig4}
\end{subfigure}
\begin{subfigure}[b]{0.24\textwidth}
                \includegraphics[width=\textwidth]{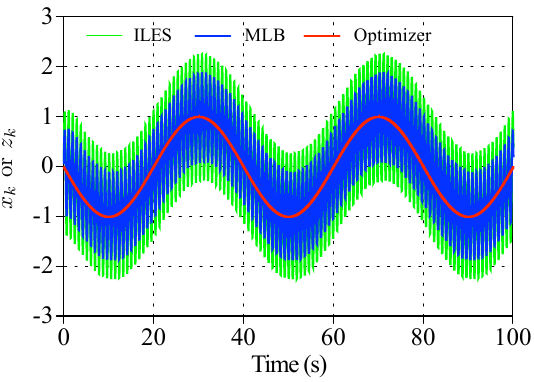}
                \subcaption{\scriptsize$50^\text{th}$ iteration with $\omega=7 {{\text{ rad/s}}}$}
                \label{fig5}
\end{subfigure}
\caption{The evolution of the original system and MLB system under frequency $\omega=7 {{\text{ rad/s}}}$.}
\end{figure*}

\begin{figure*}[htbp]
\begin{subfigure}[b]{0.24\textwidth}
                \includegraphics[width=\textwidth]{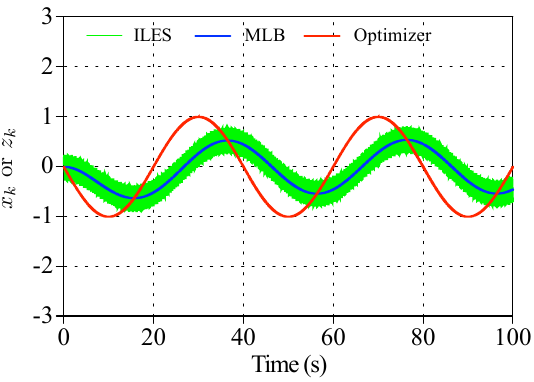}
                \subcaption{\scriptsize$1^\text{st}$ iteration with $\omega=15 {{\text{ rad/s}}}$}
                \label{fig6}
\end{subfigure}
\begin{subfigure}[b]{0.24\textwidth}
                \includegraphics[width=\textwidth]{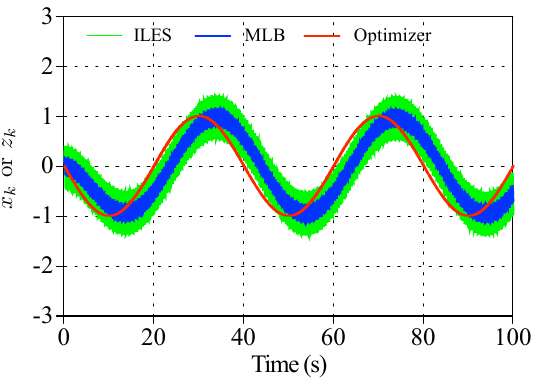}
                \subcaption{\scriptsize$2^\text{nd}$ iteration with $\omega=15 {{\text{ rad/s}}}$}
                \label{fig7}
\end{subfigure}
\begin{subfigure}[b]{0.24\textwidth}
                \includegraphics[width=\textwidth]{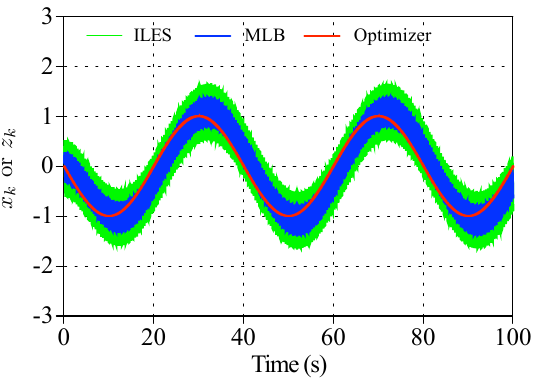}
                \subcaption{\scriptsize$3^\text{rd}$ iteration with $\omega=15 {{\text{ rad/s}}}$}
                \label{fig8}
\end{subfigure}
\begin{subfigure}[b]{0.24\textwidth}
                \includegraphics[width=\textwidth]{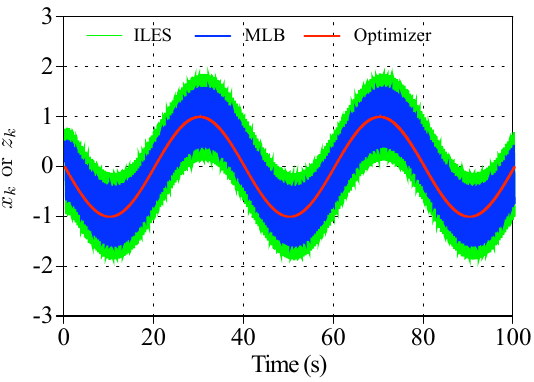}
                \subcaption{\scriptsize$50^\text{th}$ iteration with $\omega=15 {{\text{ rad/s}}}$}
                \label{fig9}
\end{subfigure}
\caption{The evolution of the original system and MLB system under frequency $\omega=15 {{\text{ rad/s}}}$.}
\end{figure*}

We study the following static map.
\begin{equation*}
y=x^2+2\sin\left(\frac{\pi}{20}t\right)x
\end{equation*}
It is evident that the minimizer trajectory is $x^*(t)=-\sin\left(\frac{\pi}{20}t\right)$. As for the iterative learning extremum seeking control law, $\alpha$ is selected as $0.1$ and the forgetting factor $\beta=0.3$. Figs. \ref{fig2}-\ref{fig5} show the evolutions of the original system and the MLB system (\ref{eqn:mlbs}) in the $1^\text{st}$, $2^\text{nd}$, $3^\text{rd}$, $50^\text{th}$ iterations under the frequency $\omega=7 {{\text{ rad/s}}}$. They indicate that the tracking performances of both the systems improve gradually, although the fluctuations are getting larger but finally bounded. Figs. \ref{fig6}-\ref{fig9} show the similar evolutions under the frequency $\omega=15 {\text{ rad/s}}$. Comparing Figs. \ref{fig2}-\ref{fig5} and Figs. \ref{fig6}-\ref{fig9}, it implies that both the tracking error (the fluctuation of the MLB system (\ref{eqn:mlbs})) and the approximation error (the gap between the original system and the MLB system (\ref{eqn:mlbs})) are smaller under a larger frequency $\omega$.


\subsection{PMFV optimization problem}\label{PMFV}
\begin{figure}[htbp]
	\begin{center}
		\includegraphics[width=5.4cm]{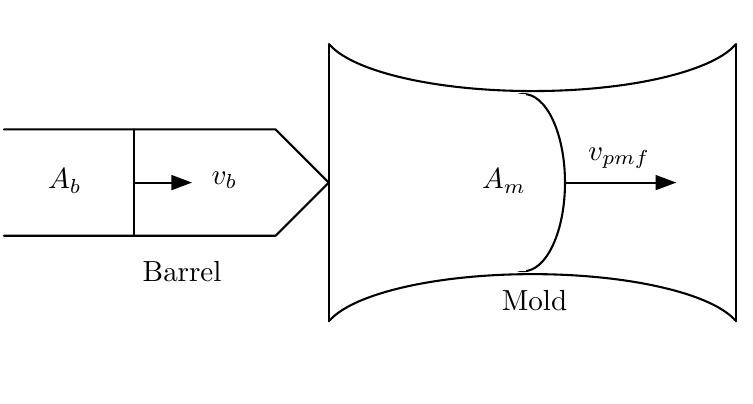}
		\caption{Relation between PMFV and IV}
		\label{fig_imm}
	\end{center}
\end{figure}

\begin{figure}[htbp]
	\begin{center}
		\includegraphics[width=0.4\textwidth]{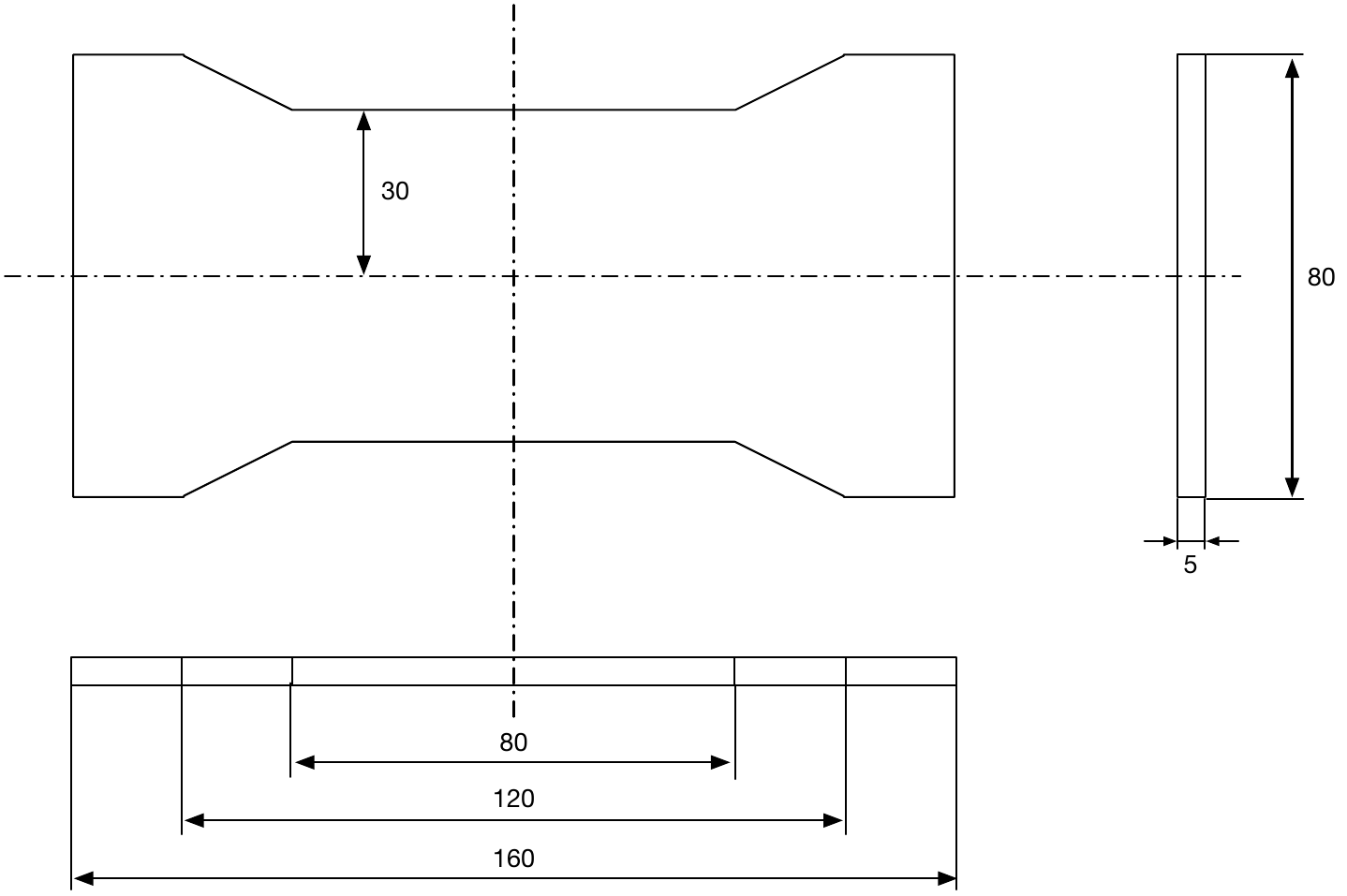}
		\caption{Geometry of the mold cavity in simulation(unit: mm)}
		\label{mold}
	\end{center}
\end{figure}

\begin{figure}[htbp]
	\begin{center}
		\includegraphics[width=0.4\textwidth]{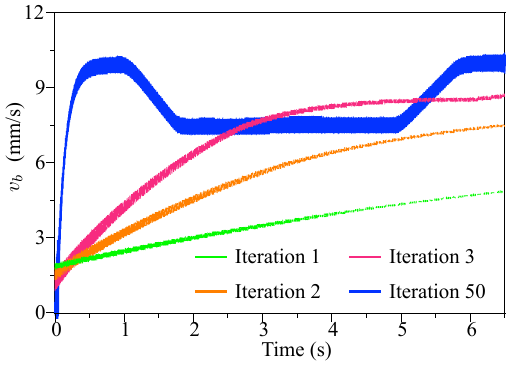}
		\caption{$v_b$ profile change versus iteration}
		\label{ES_PMFV}
	\end{center}
\end{figure}

\begin{figure}[htbp]
	\begin{center}
		\includegraphics[width=0.4\textwidth]{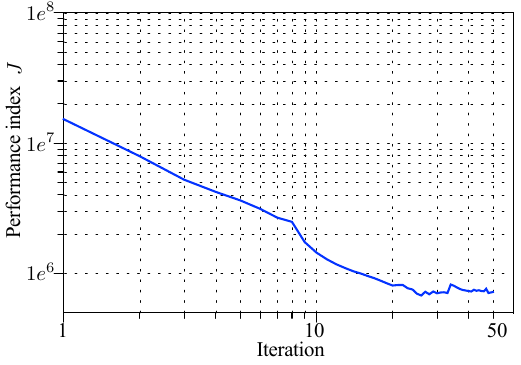}
		\caption{Performance index $J$ versus iteration in logarithm plot}
		\label{ES_PMFV_J}
	\end{center}
\end{figure}

{{Injection molding, an important polymer processing technique, transforms polymer granules into various plastic parts with high versatility and productivity. Surface quality of the products is of pivotal interest. Researches indicate that it is mainly determined by the evenness of polymer-melt-front velocity (PMFV) showed in Fig. \ref{fig_imm} and denoted as $v_{pmf}$\cite{Chen}. Fig. \ref{fig_imm} shows an abstracted injection-molding model, where $v_b$ is the injection velocity (map's input), $A_b$ the cross-section area of barrel, $A_m$ the cross-section area of the melt-front inside the mold cavity. Equally important, plastic engineers are also interested in the productivity, which can be roughly expressed as $1/v_{pmf}$. Supposing that the polymer melt is incompressible, according to Fig. \ref{fig_imm} and mass conservation, it is easy to have $v_bA_b=A_mv_{pmf}$. In practice, $A_b$ is usually known and $v_{pmf}$ can be measured by capacity transducer \cite{Chen2}. $A_m$ is a function of filling extent -- the distance between the polymer melt front and the gate. Thus, the problem is to {{steer}} $v_b$ to minimize
\begin{equation*}
J=\int_0^L\lambda_1\dot{v}_{pmf}^2+\lambda_2/v_{pmf}+\lambda_3(v_{pmf}-v_g)^2dt.
\end{equation*}
It is a composite objective function; $\lambda_1,\lambda_2$ and $\lambda_3$ are weights. From plastic engineering practice, if $v_{pmf}$ is too low, the polymer melt may get solidified before completing the filling; if $v_{pmf}$ is too high, the polymer melt may get burned due to the large shearing force exerted by the mold wall \cite{IM_handbook}. The third term $\lambda_3(v_{pmf}-v_g)^2$ is introduced to prevent these situations. $v_g$ is the guided velocity chosen by engineering experiences. Another function of $\lambda_3(v_{pmf}-v_g)^2$ is to regularize the cost function by making it convex. Heuristically, the best minimizer to $J$ is a constant $v_{pmf}^*$ such that $\dot{v}_{pmf}^2$ disappears for an appropriate set of $\lambda_1,\lambda_2,\lambda_3$. Expanding $J$ to the second order around $v_{pmf}^*$, $J\approx J(v_{pmf}^*)+ \int_0^L[\lambda_2(v_{pmf}^*)^{-3}+\lambda_3](v_{pmf}-v_{pmf}^*)^2dt$. Note that $v_b$ correlates to $v_{pmf}$ through a changing factor $A_m$; $A_m$ is a function of displacement, i.e., $\phi(v_{pmf}^*t)$. Thus, $J$ admits a time-varying minimizer $v_b^*=Kv_{pmf}^*\phi(v_{pmf}^*t)$ for some constant $K$.
The mold for the numerical study is shown in Fig. \ref{mold} \cite{mold}. $\lambda_1=5\times 10^{-6},\lambda_2=30,\lambda_3=1,\beta=0.09,\omega=1000\text{ rad/s}, \alpha=0.01, v_g=25 \text{ mm/s}, Ab=1000 \text{ mm}^2$. Fig. \ref{ES_PMFV} shows the barrel velocity $v_b$ profile change versus iteration. Fig. \ref{ES_PMFV_J} demonstrates the performance index $J$ decreases as the iteration goes. It also suggests that the convergence rate of our algorithm is exponentially fast.}}

\section{Conclusion \& Outlook}
This paper has proposed an iterative learning extremum seeking approach to solve {{{the optimization problem for repetitive time-varying mapping}}. A modified Lie bracket system has been introduced to study the behavior of the ILES system. 
It has shown that the MLB system is an online integral-type ILC control law with the bounded approximation error. The convergence of the corresponding ILC control law has been analyzed. Based on that, the convergence of the proposed ILES to a set has been shown. The size of the set is reducible by tuning the frequency of the dither signal. {{The distance from set's center ($y_\infty$) to the origin ($0$)}} is also tunable by some appropriate forgetting factors $\beta$. 
In the future, it is quite interesting to investigate how to extend the method to cover a general function $F(x,t)$. It is worthy to study the dynamic mapping situation as well.

% or


\appendices
\section*{Appendix}
\subsection*{Proof of Theorem \ref{thm-4-1}}\label{pf_thm_4_1}
\begin{IEEEproof}
The proof is similar to Theorem B.1 in \cite{Khalil}. Arbitrarily select $x_1\in S$ and generate a sequence $\{x_k\}$ according to the formula $x_{k+1}=T_k(x_k)$. Every $x_k\in S$, since $T_k\in\mathsf{B}(S)$. First, we will show $\{x_k\}$ is a Cauchy sequence. It follows from the definition of $S$ in (\ref{eqn:set_S}) that there is a constant $D>0$ such that $\|x\|\le D$ for all $x\in S$. Additionally, that $\{T_k\}$ is a Cauchy sequence follows, since $T_k$ converges uniformly. Then, for an arbitrary $\epsilon>0$, there exists a $k_\epsilon$ such that $\|T_m-T_n\|\le \epsilon/D$ for any $m,n\ge k_\epsilon$. For $k-1>k_\epsilon$, we have
\begin{equation*}
\begin{aligned}
\|x_{k+1}-x_k\|=&\|T_k(x_k)-T_{k-1}(x_{k-1})\|\\
\le&\|T_k(x_k)-T_k(x_{k-1})\|\\
&+\|T_k(x_{k-1})-T_{k-1}(x_{k-1})\|\\
\le&\rho\|x_k-x_{k-1}\|+\|T_k-T_{k-1}\|\|x_{k-1}\|\\
\le&\rho\|x_k-x_{k-1}\|+\epsilon
\end{aligned}
\end{equation*}
It follows that
\begin{equation*}
\begin{aligned}
\|x_{k+r}-x_k\|\le&\sum_{i=0}^{r-1}\|x_{k+i+1}-x_{k+i}\|\\
\le&\sum_{i=0}^{r-1}\left[\rho^{i+1}\|x_{k}-x_{k-1}\|+\rho^i\epsilon\right]\\
\le&\frac{\rho}{1-\rho}\|x_k-x_{k-1}\|+\frac{\epsilon}{1-\rho}\\
\le&\frac{\rho}{1-\rho}(\rho\|x_{k-1}-x_{k-2}\|+\epsilon)+\frac{\epsilon}{1-\rho}\\
\le&\frac{\rho}{1-\rho}\left(\rho^{k-k_\epsilon-1}\|x_{k_\epsilon+1}-x_{k_\epsilon}\|+\sum_{i=0}^{k-k_\epsilon-2}\rho^{i}\epsilon\right)\\
&+\frac{\epsilon}{1-\rho}\\
\le&\frac{2\rho^{k-k_\epsilon}D}{1-\rho}+\frac{\rho\epsilon}{(1-\rho)^2}+\frac{\epsilon}{1-\rho}\\
\le&\frac{2\rho^{k-k_\epsilon}D}{1-\rho}+\frac{\epsilon}{(1-\rho)^2}
\end{aligned}
\end{equation*}
Since $\epsilon$ is chosen arbitrarily, the right hand side will go to $0$ as $k\to\infty$. Hence, $\{x_k\}$ is a Cauchy sequence. From that $\mathcal{X}$ is complete, $x_k\to x^\star\in\mathcal{X}$ as $k\to\infty$. Furthermore, $S$ is closed; it follows that $x^\star\in S$.

%{{The next steps are to show that $x^\star=T_\infty(x^\star)$ and the uniqueness {{of the limit solution $x^\star$}}. Their proofs follow the ideas from \cite{Khalil}, p. 656.}}
The second step is to prove $x^\star=T_\infty(x^\star)$. For any $x_{k+1}=T_k(x_k)$, it can be obtained that
\begin{equation*}
\begin{aligned}
\|x^\star-T_\infty(x^\star)\|\le&\|x^\star-x_{k+1}\|+\|x_{k+1}-T_\infty(x^\star)\|\\
\le&\|x^\star-x_{k+1}\|+\|x_{k+1}-T_\infty(x_k)\|\\
&+\|T_\infty(x_k)-T_\infty(x^\star)\|\\
\le&\|x^\star-x_{k+1}\|+\|T_k-T_\infty\|\|x_k\|\\
&+\rho\|x_k-x^\star\|
\end{aligned}
\end{equation*}
It is apparent that the right hand side can be made arbitrarily small by selecting a sufficiently large $k$. Thereafter, $x^\star=T_\infty(x^\star)$.

Finally, we will show the uniqueness. Suppose that there is another fix point $y^\star$ satisfying $y^\star=T_\infty(y^\star)$. Then, we have
\begin{equation*}
\|x^\star-y^\star\|\le\|T_\infty(x^\star)-T_\infty(y^\star)\|\le\rho\|x^\star-y^\star\|
\end{equation*}
Since $\rho$ is strictly less than $1$, it is a contradiction; then $x^\star=y^\star$. This completes the whole proof.
\end{IEEEproof}

\subsection*{Proof of Lemma \ref{lem1}}\label{pf1}
\begin{IEEEproof}
From the definitions of $T_k$ and $\lambda$-norm, we have
\begin{equation*}
\begin{aligned}
&\|T_k(x)-T_k(y)\|_\lambda\\
=&\max_{t\in[0,L]}\text{e}^{-\lambda t}(1-\beta)\\
&\left\|(x-y)-\text{e}^{-\Gamma_kt}\int_0^t\text{e}^{\Gamma_ks}\Gamma_k(x-y)ds\right\|_\infty
\end{aligned}
\end{equation*}
It follows from the triangle inequality of norm and $\max\{a+b\}\le\max\{a\}+\max\{b\}$ that
\begin{equation*}
\begin{aligned}
&\|T_k(x)-T_k(y)\|_\lambda\\
\le&\max_{t\in[0,L]}\text{e}^{-\lambda t}(1-\beta)\|x-y\|_\infty\\
&+ \max_{t\in[0,L]}\text{e}^{-\lambda t}(1-\beta)\\
&\left\|\text{e}^{-\Gamma_k t}\int_0^t \text{e}^{\Gamma_k s}\Gamma_k(x-y)ds\right\|_\infty
\end{aligned}
\end{equation*}
The first term on the right hand side is exactly $(1-\beta)\|x-y\|_\lambda$ according to the definition of $\lambda$-norm. For the second term, we do the following operation.
\begin{equation*}
\begin{aligned}
&\|T_k(x)-T_k(y)\|_\lambda\le(1-\beta)\|x-y\|_\lambda+ \max_{t\in[0,L]}\text{e}^{-\lambda t}(1-\beta)\\
&\left\|\text{e}^{-\Gamma_kt}\int_0^t \text{e}^{(\Gamma_k +\lambda I)s}\text{e}^{-\lambda s}\Gamma_k(x-y)ds\right\|_\infty
\end{aligned}
\end{equation*}
According to mean value theorem, it can be obtained that
\begin{equation*}
\begin{aligned}
&\left\|\text{e}^{-\Gamma_kt}\int_0^t \text{e}^{(\Gamma_k +\lambda I)s}\text{e}^{-\lambda s}\Gamma_k(x-y)ds\right\|_\infty\\
=&\left\|\text{e}^{-\Gamma_kt}\int_0^t \text{e}^{(\Gamma_k +\lambda I)s}\Gamma_kds\left(\text{e}^{-\lambda \xi}(x(\xi)-y(\xi))\right)\right\|_\infty,\\
&[\xi\in(0,t)]\\
\le&\text{e}^{-\lambda\xi}\|x(\xi)-y(\xi)\|_\infty\left\|\text{e}^{-\Gamma_kt}\int_0^t\text{e}^{(\Gamma_k+\lambda I)s}\Gamma_kds\right\|_\infty,\\
&[\xi\in(0,t)]\\
\le&\max_{s\in[0,t]}\left\{\text{e}^{-\lambda s}\|x-y\|_\infty\right\}\left\|\text{e}^{-\Gamma_kt}\int_0^t\text{e}^{(\Gamma_k+\lambda I)s}\Gamma_kds\right\|_\infty
\end{aligned}
\end{equation*}
The first inequality is followed from the definition of induced matrix norm. It is also noted that
\begin{equation*}
\max_{s\in[0,t]}\left\{\text{e}^{-\lambda s}\|x-y\|_\infty\right\}\le\|x-y\|_\lambda
\end{equation*}
Hence, 
\begin{equation*}
\begin{aligned}
&\|T_k(x)-T_k(y)\|_\lambda\\
\le&(1-\beta)\|x-y\|_\lambda\left(1+\max_{t\in[0,L]}\text{e}^{-\lambda t}\right.\\
&\left.\left\|\text{e}^{-\Gamma_kt}\int_0^t\text{e}^{(\Gamma_k+\lambda I)s}\Gamma_kds\right\|_\infty\right)\\
\le&(1-\beta)\|x-y\|_\lambda\left(1+\max_{t\in[0,L]}\left\|\int_0^t\text{e}^{(\Gamma_k+\lambda I)(s-t)}\Gamma_kds\right\|_\infty\right)\\
\le&(1-\beta)\|x-y\|_\lambda\left[1+\max_{t\in[0,L]}\left\|\left(I-\text{e}^{-(\Gamma_k+\lambda I)t}\right)\Gamma_k\right.\right.\\
&(\Gamma_k+\lambda I)^{-1}\Big\|_\infty\bigg]
\end{aligned}
\end{equation*}
For a matrix $A\in\mathbb{M}_n$, $\frac{1}{\sqrt{n}}\|A\|_2\le\|A\|_\infty\le\sqrt{n}\|A\|_2$, which is obtained from 
the norm equivalence theorem. Then, it follows that
\begin{equation*}
\begin{aligned}
&\left\|\left(I-\text{e}^{-(\Gamma_k+\lambda I)t}\right)\Gamma_k(\Gamma_k+\lambda I)^{-1}\right\|_\infty\\
\le&\sqrt{n}\left\|\left(I-\text{e}^{-(\Gamma_k+\lambda I)t}\right)\Gamma_k(\Gamma_k+\lambda I)^{-1}\right\|_2\\
\le&\sqrt{n}\left\|\Gamma_k(\Gamma_k+\lambda I)^{-1}\right\|_2
\end{aligned}
\end{equation*}
Thereafter, we have
\begin{equation*}
\|T_k(x)-T_k(y)\|_\lambda\le(1-\beta)\|x-y\|_\lambda(1+\sqrt{n}\|\Gamma_k(\Gamma_k+\lambda I)^{-1}\|_2)
\end{equation*}
Note that $\|\Gamma_k(\Gamma_k+\lambda I)^{-1}\|_2=\|I-(\Gamma_k+\lambda I)^{-1}\|_2$, it is a non increasing function of $\Gamma_k$. From the assumption, it is known that $\Gamma_k\ge \alpha\delta/2 I$. Denoting $\rho=(1-\beta)(1+\sqrt{n}\alpha\delta/(\alpha\delta+2\lambda))$, it is easy to see there always exists a $\rho<1$ when $\lambda\in(\lambda_0,+\infty)$ for arbitrary $\beta\in(0,1)$ with
\begin{equation*}
\lambda_0=\max\left\{0,\frac{\alpha\delta[\sqrt{n}(1-\beta)-\beta]}{2\beta}\right\}
\end{equation*}
This completes the proof.
%Similar conclusions can be drawn on arbitrary $\lambda$.
\end{IEEEproof}

\subsection*{Proof of Lemma \ref{lem_ilc1}}\label{pf_ilc1}
\begin{IEEEproof}
We are going to show the statement by inductive arguments. 

When $k=1$, from (\ref{eqn:exp_sol}), the explicit solution of $z_1$ is
\begin{equation*}
z_1=-\int_0^t\text{e}^{-\Gamma_1(t-s)}\Gamma_1x^*ds
\end{equation*}
Since all the eigenvalues associated with $\Gamma_1$ locates in the right complex plane and $x^*$ is continuous over the interval $[0,L]$, $z_1$ is well defined and bounded over $[0,L]$. Hence, $\|y_1\|_\lambda$ is bounded. Furthermore, it is evident that $\|y_1\|_\lambda\le D_2$, equivalently, $y_1\in S$.

Assume that $y_{k-1}\in S$ or $\|y_{k-1}\|_\lambda\le D_2$ holds for the {{$(k-1)$}}-th iteration. It remains to show that $y_{k}\in S$.

From (\ref{eqn:yk3}) and Lemma \ref{lem1}, we have
\begin{equation*}
\begin{aligned}
\|y_k\|_\lambda=&\|G_k(y_{k-1})\|_\lambda \\
\le& \left\|T_k(y_{k-1})-T_k(0)+T_k\left(-\frac{\beta}{1-\beta}x^*\right)-T_k(0)\right\|_\lambda\\
\le&\|T_k(y_{k-1})-T_k(0)\|_\lambda\\
&+\left\|T_k\left(-\frac{\beta}{1-\beta}x^*\right)-T_k(0)\right\|_\lambda\\
\le & \rho\|y_{k-1}\|_\lambda+\rho\left\|\frac{\beta}{1-\beta}x^*\right\|_\lambda\\
\le & \rho D_2+\rho\left\|\frac{\beta}{1-\beta}x^*\right\|_\lambda
\end{aligned}
\end{equation*}
From the statement, we have $D_2\ge\frac{\beta\rho}{(1-\beta)(1-\rho)}\|x^*\|_\lambda$. Thus,
\begin{equation*}
\|y_k\|_\lambda\le\rho D_2+\rho\left\|\frac{\beta}{1-\beta}x^*\right\|_\lambda\le \rho D_2+(1-\rho)D_2=D_2
\end{equation*}
Hence, $G_k(y_{k-1})\in S$, and this completes the proof.
\end{IEEEproof}

\subsection*{Proof of Theorem \ref{thm_ilc1}}\label{pf_thm_ilc1}
{{
\begin{IEEEproof}
It is {{noted that (\ref{eqn:y_inf})}} is equivalent to 
\begin{equation*}
y_\infty=G_\infty(y_\infty)
\end{equation*}
Hence, we only need to check whether $G_k$ satisfies C1-C3 or not. Since $\mathcal{X}$ is equipped with $\lambda$-norm, we can define the mapping norm by the $\lambda$-norm as 
\begin{equation*}
\|G_k\|=\sup_{\|x\|_\lambda=1}\|G_k(x)\|_\lambda
\end{equation*}
Expanding $G_k(x)$, we have
\begin{equation*}
\|G_k\|\le\sup_{\|x\|_\lambda=1}\|T_k(x)\|_\lambda+\left\|T_k\left(-\frac{\beta}{1-\beta}x^*\right)\right\|_\lambda
\end{equation*}
On the other hand, from Lemma \ref{lem1}, it is known that $\|T_k(x)\|_\lambda=\|T_k(x)-T_k(0)\|_\lambda\le\rho\|x\|_\lambda$. It immediately follows that
\begin{equation*}
\|G_k\|\le\rho+\left\|T_k\left(-\frac{\beta}{1-\beta}x^*\right)\right\|_\lambda<+\infty
\end{equation*}
Combining with Lemma \ref{lem_ilc1}, $G_k\in \mathsf{B}(S)$ for arbitrary $k$.
Also noted from Lemma \ref{lem1}, the mapping sequence $\{G_k\}$ satisfies C1 and C2; it remains to show that $G_k$ converges to $G_\infty$ uniformly. It suffices to show $T_k$ converges to $T_\infty$ uniformly. By definition, we have
\begin{equation*}
\|T_k-T_\infty\|=\sup_{\|x\|_\lambda=1}\|T_k(x)-T_\infty(x)\|_\lambda
\end{equation*}
From the definition of $\lambda$-norm and (\ref{eqn:yk1}), we further get that
\begin{equation*}
\begin{aligned}
&\|T_k-T_\infty\|\le(1-\beta)\\
&\times\sup_{\substack{t\in[0,L]\\\|x\|_\lambda=1}}\text{e}^{-\lambda t}\left\|\int_0^t\left(\text{e}^{-\Gamma_k(t-s)}\Gamma_k-\text{e}^{-\Gamma_\infty(t-s)}\Gamma_\infty\right)xds\right\|_\infty
\end{aligned}
\end{equation*}
By the mean value theorem, it is easy to derive that
{\small
\begin{equation*}
\begin{aligned}
&\|T_k-T_\infty\|\\
\le&(1-\beta)\sup_{\substack{t\in[0,L]\\\|x\|_\lambda=1}}\text{e}^{-\lambda t}\left\|\int_0^t\left(\text{e}^{-\Gamma_k(t-s)}\Gamma_k-\text{e}^{-\Gamma_\infty(t-s)}\Gamma_\infty\right)ds\right\|_\infty\\
&\times\|x(\xi)\|_\infty, (\xi\in[0,L])\\
\le&(1-\beta)\sup_{t\in[0,L]}\left\|\int_0^t\left(\text{e}^{-\Gamma_k(t-s)}\Gamma_k-\text{e}^{-\Gamma_\infty(t-s)}\Gamma_\infty\right)ds\right\|_\infty\\
\le&(1-\beta)\sup_{t\in[0,L]}\left\|\text{e}^{-\Gamma_k t}-\text{e}^{-\Gamma_\infty t}\right\|_\infty
\end{aligned}
\end{equation*}}The term $\exp(-\Gamma_k t)$ can approach $\exp(-\Gamma_\infty t)$ arbitrarily small as $k\to\infty$. Therefore, it can be concluded that
\begin{equation*}
\|T_k-T_\infty\|\to 0 \text{  as  } k\to\infty
\end{equation*}
which is the uniform convergence. Then, apply the $k$-dependent contraction mapping, we can conclude the result.
\end{IEEEproof}
}}

\subsection*{Proof of Thoerem \ref{thm_ilc2}}\label{pf_thm_ilc2}
\begin{IEEEproof}
From (\ref{eqn:mlb0z}), within this proof, we are equivalently studying
\begin{equation}\label{eqn:thm4-2-1}
\dot{y}_k=\dot{y}_{k-1}-\Gamma y_k
\end{equation}
Define an index as follows.
\begin{equation}\label{eqn:thm4-2-2}
J_k=\int_0^L\text{e}^{-\lambda t}\dot{y}_k^T\dot{y}_kdt
\end{equation}
where $\lambda>0$. Rewriting (\ref{eqn:thm4-2-1}) as $\dot{y}_{k-1}=\dot{y}_k+\Gamma y_k$, we insert it into (\ref{eqn:thm4-2-2}) and compare the difference of (\ref{eqn:thm4-2-2}) between $k$ and $k-1$.
\begin{equation*}
\begin{aligned}
J_k-J_{k-1}=&\int_0^L\text{e}^{-\lambda t}\left[\dot{y}_k^T\dot{y}_k-(\dot{y}_k+\Gamma y_k)^T(\dot{y}_k+\Gamma y_k)\right]dt\\
=&-\int_0^L\text{e}^{-\lambda t}y_k^T\Gamma^T\Gamma y_kdt-2\int_0^L\text{e}^{-\lambda t}y_k\Gamma\dot{y}_kdt
\end{aligned}
\end{equation*}
As for the second term in the right hand side, we integrate by parts and derive that
\begin{equation*}
\begin{aligned}
2\int_0^L\text{e}^{-\lambda t}y_k^T\Gamma\dot{y}_kdt=&\text{e}^{-\lambda t}y_k^T\Gamma y_k\bigg|_{t=0}^{t=L}+\lambda\int_0^L\text{e}^{-\lambda t}y_k^T\Gamma y_kdt
\end{aligned}
\end{equation*}
Thus, we have
\begin{equation*}
\begin{aligned}
J_k-J_{k-1}=&-\int_0^L\text{e}^{-\lambda t}y_k^T\left(\Gamma^T\Gamma+\lambda\Gamma\right)y_kdt\\
&-\text{e}^{-\lambda t}y_k^T(L)\Gamma y_k(L)\\
\le&-\rho_{\min}\int_0^L\text{e}^{-\lambda t}y_k^Ty_kdt
\end{aligned}
\end{equation*}
Let $\rho_{\min}$ be the smallest eigenvalue of the matrix $\Gamma^T\Gamma+\lambda\Gamma$. It is apparent that $\rho_{\min}>0$ since $\Gamma$ is positive definite and $\lambda>0$. Thus, $\{J_k\}$ is a non-increasing real number sequence, and $J_k$ is bounded below by $0$, $J_k$ converges. It follows that
\begin{equation*}
\lim_{k\to\infty}\left(\int_0^Ly_k^Ty_kdt\right)^{\frac{1}{2}}=0
\end{equation*}
It is the $L_2$-norm of $y_k$. Thus, we can claim that $y_k$ converges to $0$ almost everywhere.
\end{IEEEproof}

\subsection*{Proof of Lemma \ref{lem2}}\label{pf2}
\begin{IEEEproof}
The proof is adapted from Theorem 1 in \cite{Durr}. The major differences are: first, \cite{Durr} dealing with more general form -- input affine system, we only focus on sinusoid input signal; second, \cite{Durr} suitable for infinite time horizon and single iteration, we are handling the case of finite time horizon and multiple iterations.

\begin{figure}[htbp]
                \includegraphics[width=0.45\textwidth]{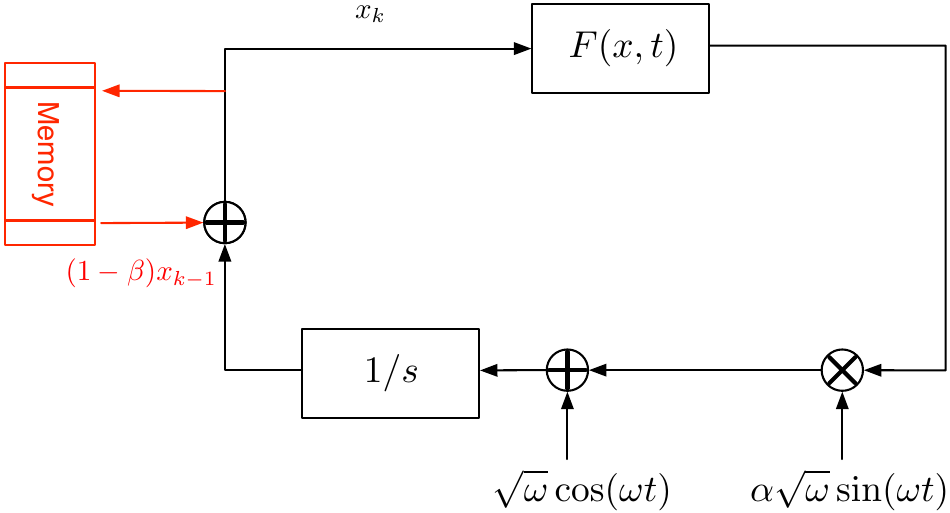}
                \caption{The illustration of the proof idea}
                \label{ILESpf}
\end{figure}

The basic idea to show the lemma is illustrated in Fig. \ref{ILESpf}. $x_{k< k_0}$ are within $S$ and $z_{k_0}$ is within $\mathcal{K}$. We construct a tube with radius $E$ along the trajectory of $z_{k_0}$. If $E\le\text{dist}(\mathcal{K},S)$ and $x_{k_0}$ is within the tube, it can conclude that $x_{k_0}\in S$. Moreover, supposing the $x_{k_0}$ leaves the tube at arbitrary time $t_E$, if we can show $t_E$ is not the time that $x_{k_0}$ leaves the tube, then we can claim that $x_{k_0}$ will never leave the tube.

Evaluating (\ref{eqn:iles1}), (\ref{eqn:mlbs}) at $k_0$, Subtracting (\ref{eqn:mlbs}) from (\ref{eqn:iles1}) and integrating on both sides, we have
\begin{equation}\label{eqn:xz}
\begin{aligned}
x_{k_0}-z_{k_0}=&-\alpha\int_0^tF(x_{k_0},t)\sqrt{\omega}\sin(\omega s)ds\\
&+\int_0^t\sqrt{\omega}\cos(\omega s)ds\\
&+\frac{\alpha\gamma_{k_0}}{2}\int_0^t\nabla_zF(z_{k_0},s)ds
\end{aligned}
\end{equation}
Let $R_a$ be the first term on the right hand side of (\ref{eqn:xz}); taking integration on $R_a$ by parts, we obtain that
\begin{equation*}
\begin{aligned}
R_a=&\frac{\alpha}{\sqrt{\omega}}\left.\left[F(x_{k_0},s)\cos(\omega s)\right]\right|_0^t-\frac{\alpha}{\sqrt{\omega}}\int_0^t\cos(\omega s)\dot{F}(x_{k_0},s)ds\\
=&\frac{\alpha}{\sqrt{\omega}}\left.\left[F(x_{k_0},s)\cos(\omega s)\right]\right|_0^t\\
&-\frac{\alpha}{\sqrt{\omega}}\int_0^t\cos(\omega s)\frac{\partial F(x_{k_0},s)}{\partial t}ds\\
&-\frac{\alpha}{\sqrt{\omega}}\int_0^t\cos(\omega s)\nabla_xF(x_{k_0},s)\dot{x}_{k_0}ds
\end{aligned}
\end{equation*}
Denote the $1^{\text{st}}$ and $2^{\text{nd}}$ terms in the right hand side of the above equation as $R_1$ and $R_2$ respectively. Iterating (\ref{eqn:iles1}) to the first iteration, it follows that
\begin{equation}\label{eqn:iter}
\dot{x}_{k_0}=-\sum_{i=1}^{k_0}(1-\beta)^{k_0-i}\alpha F(x_i,t)\sqrt{\omega}\sin(\omega t)+\gamma_{k_0}\sqrt{\omega}\cos(\omega t)
\end{equation}
Inserting (\ref{eqn:iter}) into $R_a$, it becomes
\begin{equation*}
\begin{aligned}
R_a=&R_1+R_2+R_3-\alpha\gamma_{k_0}\int_0^t\cos^2(\omega s)\nabla_x F(x_{k_0},s)ds\\
=&R_1+R_2+R_3+R_4-\frac{\alpha\gamma_{k_0}}{2}\int_0^t\nabla_xF(x_{k_0},s)ds
\end{aligned}
\end{equation*}
Here $R_3\triangleq\frac{\alpha^2}{2}\sum_{i=1}^{k_0}(1-\beta)^{k_0-i}\int_0^t\sin(2\omega s)\nabla_xF(x_{k_0},s)F(x_i,s)ds$ and $R_4\triangleq-\frac{\alpha\gamma_{k_0}}{2}\int_0^t\cos(2\omega s)\nabla_xF(x_{k_0},s)ds$. The second equality stems from the identity $\cos(2\omega s)=2\cos^2(\omega s)-1$. Hence, $x_{k_0}-z_{k_0}$ becomes
\begin{equation}\label{eqn:xz1}
x_{k_0}-z_{k_0}=\sum_{i=1}^5R_i-\frac{\alpha\gamma_{k_0}}{2}\int_0^t[\nabla_xF(x_{k_0},s)-\nabla_zF(z_{k_0},s)]ds
\end{equation}
where $R_5=\sin(\omega t)/{\sqrt{\omega}}$. It should be noted that (\ref{eqn:xz1}) holds universally for any $t$. 

{{Now}} we are going to show $x_{k_0}$ will remain in $\mathcal{S}$ over $[0,L]$ by contradiction. Let $E=\text{dist}(\mathcal{K},\mathcal{S})$. Assume that there is a time instant $t_E\in(0,L)$ such that $\|x_{k_0}(t)-z_{k_0}(t)\|_\infty<E$ for any $t\in[0,t_E)$ and $\|x_{k_0}(t_E)-z_{k_0}(t_E)\|_\infty=E$. It means that $t_E$ is the first time when $x_{k_0}$ leaves a tube with $z_{k_0}$ as center and radius $E$. Since $\mathcal{S}$ is a compact set, and $f^*(t),x^*(t)\in C_1[0,L]$, by the definition of $F(x,s)$, for any $x\in\mathcal{S}$, there always exist $\epsilon_1,\epsilon_2,\epsilon_3$ such that
\begin{equation*}
\|F(x,s)\|_\infty\le \epsilon_1, \|\nabla_xF(x,s)\|_\infty\le\epsilon_2,\left\|\frac{\partial F(x,s)}{\partial t}\right\|_\infty\le\epsilon_3.
\end{equation*}
Note that these bounds are uniformly bounded in $t$. 
Thus, we have
\begin{equation*}
\begin{aligned}
\|R_1\|_\infty=&\left\|\frac{\alpha}{\sqrt{\omega}}\left.\left[F(x_k,s)\cos(\omega s)\right]\right|_0^t\right\|_\infty\le\frac{2\alpha\epsilon_1}{\sqrt{\omega}}\\
\|R_2\|_\infty=&\left\|\frac{\alpha}{\sqrt{\omega}}\int_0^t\cos(\omega s)\frac{\partial F(x_k,s)}{\partial t}ds\right\|_\infty\le\frac{\alpha\epsilon_3}{\omega^{3/2}}\\
\|R_3\|_\infty=&\left\|\frac{\alpha^2}{2}\sum_{i=1}^{k_0}(1-\beta)^{k_0-i}\int_0^t\sin(2\omega s)\nabla_xF(x_k,s)\right.\\
&F(x_i,s)ds\bigg\|_\infty\le\frac{\alpha^2\epsilon_1\epsilon_2}{2\omega}\sum_{i=1}^{k_0}(1-\beta)^{k_0-i}\\
&<\frac{\alpha^2\epsilon_1\epsilon_2}{2\beta\omega}\\
\|R_4\|_\infty=&\left\|\frac{\alpha\gamma_{k_0}}{2}\int_0^t\cos(2\omega s)\nabla_xF(x_{k_0},s)ds\right\|_\infty\le\frac{\alpha\gamma_{k_0}\epsilon_2}{4\omega}\\
&\le\frac{\alpha\epsilon_2}{4\beta\omega}\\
\|R_5\|_\infty=&\left\|\frac{\sin(\omega s)}{\sqrt{\omega}}\right\|_\infty\le\frac{1}{\sqrt{\omega}}
\end{aligned}
\end{equation*}
So there must exist a positive number $M$ such that $\sum_{i=1}^5R_i\le\frac{M}{\sqrt{\omega}}$ and $M$ is independent of $k$.

Since $F(x,t)$ are twice continuous on $\mathcal{S}\times[0,t_E]$, there exist a positive number $K$ such that
$\|\nabla_xF(x_{k_0},s)-\nabla_zF(z_{k_0},s)\|_\infty\le K\|x_{k_0}-z_{k_0}\|_\infty$ according to Lemma 3.2 (pp. 90, \cite{Khalil}). Thereafter, for $t\in[0,t_E]$
\begin{equation*}
\|x_{k_0}-z_{k_0}\|_\infty\le\frac{M}{\sqrt{\omega}}+\frac{\alpha K}{2\beta}\int_0^t\|x_{k_0}-z_{k_0}\|_\infty ds
\end{equation*}
From Gronwall-Bellman inequality (pp. 651, \cite{Khalil}), we have
\begin{equation}\label{eqn:xzk}
\|x_{k_0}-z_{k_0}\|_\infty\le \frac{M}{\sqrt{\omega}}\text{e}^{\frac{\alpha K}{2\beta}t}, t\in[0,t_E]
\end{equation}
For any $\omega_0 \in(M^2\text{e}^{\frac{\alpha Kt_E}{\beta}}/{E^2},\infty)$, where $\omega_0$ is independent of $k$ and $t$, $\|x_{k_0}(t_E)-z_{k_0}(t_E)\|_\infty < E$, which contradicts. Thus, $x_{k_0}$ will remain in the tube of $z_{k_0}$. Since $z_{k_0}$ stays in $\mathcal{K}$, $x_{k_0}$ will stay in $\mathcal{S}$. (\ref{eqn:xzk}) suggests that $\|x_k-z_k\|_\lambda$ can be small enough if $\omega$ is large enough. This completes the proof.
\end{IEEEproof}

\subsection*{Proof of Theorem \ref{thm}}\label{pf_thm}
\begin{IEEEproof}
We are going to show the theorem in an inductive way.

First, at $k=1$, both the extremum seeking system {{(\ref{eqn:iles1})}} and the MLB system (\ref{eqn:mlbs}) are in the standard form. From (\ref{eqn:exp_sol}), the explicit solution of the MLB system (\ref{eqn:mlbs}) is
\begin{equation*}
z_1=\int_0^t\text{e}^{-\Gamma_1(t-s)}\Gamma_1x^*ds
\end{equation*}
Because all the eigenvalues of $-\Gamma_1$ lie in the open left complex plane and $x^*$ is contained in a compact set, $z_1$ is well defined {{for initial condition $z_1(0)=0$}} in a compact set over the time interval $[0,L]$. According to the definition of $D_2$, it is known that $\|y_1\|_\lambda<D_2$. Thus, by Lie bracket theorem \cite{Durr}, it is known that the distance between $x_1$ and $z_1$ can be arbitrarily small provided that the frequency $\omega$ is sufficiently large. Hence, there exists a $\omega_1$ such that for every $\omega\in(\omega_1,\infty)$, we have $\|x_1-z_1\|_\lambda<D_1$.

Now, we assume that $\|y_i\|_\lambda<D_2$ and $\|x_i-z_i\|_\lambda<D_1$ for any $i=1,2,\ldots, k-1$ over the entire time interval. Then, we will show that these two relations are still valid for $i=k$. From Lemma \ref{lem2}, it is easy to conclude that there exists a $\omega_2$ such that for every $\omega\in(\omega_2,+\infty)$, $\|x_k-z_k\|_\lambda<D_1$. We only need to show $\|y_k\|_\lambda<D_2$.

By the linearity of $T_k$, we have
\begin{equation*}
\begin{aligned}
\|y_k\|_\lambda=&\bigg\|H_k(y_{k-1})-H_k(0)\\
&\left.+T_k\left(x_{k-1}-z_{k-1}-\frac{\beta}{1-\beta}x^*\right)\right\|_\lambda\\
\le&\left\|T_k\left(x_{k-1}-z_{k-1}-\frac{\beta}{1-\beta}x^*\right)-T_k(0)\right\|_\lambda\\
&+\|H_k(y_{k-1})-H_k(0)\|_\lambda\\
\le & \rho\|y_{k-1}\|_\lambda+\rho\left\|x_{k-1}-z_{k-1}-\frac{\beta}{1-\beta}x^*\right\|_\lambda\\
\le & \rho D_2+\rho \left(D_1+\frac{\beta}{1-\beta}\|x^*\|_\lambda\right)\le D_2
\end{aligned}
\end{equation*}
Therefore, by selecting $\omega_0=\max\{\omega_1,\omega_2\}$, we have $\|y_k\|_\lambda<D_2$ for all $k$. This completes the proof.
\end{IEEEproof}

\subsection*{Proof of Theorem \ref{thm_4_2}}\label{pf_4_2}
\begin{IEEEproof}
From Lemma \ref{lem2} and Theorem \ref{thm}, we know that there exists a $\omega_0$ for every $\omega\in(\omega_0,+\infty)$ such that we can ensure that $\|x_k-z_k\|_\lambda\le D_1$ and $y_k\in S$ for any $k$. %Additionally, we have $D_1=\mathcal{O}\left(\frac{1}{\sqrt{\omega}}\right)$.

By similar arguments in the proof of Theorem \ref{thm_ilc1}, from the definition of limit, for an arbitrarily chosen $\epsilon>0$, selecting an $\epsilon'>0$ satisfying $2\epsilon'/(1-\rho)<\epsilon$, 
\begin{equation*}
\frac{2\epsilon'D_2}{1-\rho}<\epsilon
\end{equation*}
there exists a $k_{\epsilon}$ such that $\|G_k-G_\infty\|<\epsilon'$ can be guaranteed as long as $k\ge k_{\epsilon}$.

From (\ref{eqn:set_y}), $\mathcal{Y}$ is a $\lambda$-norm ball with $y_\infty$ as its center and $D_y$ as its radius. Hence, we have
\begin{equation*}
\text{dist}(y_k,\mathcal{Y})=\max\left\{\|y_k-y_\infty\|_\lambda-D_y,0\right\}
\end{equation*}
If $\text{dist}(y_k,\mathcal{Y})=0$, it means that $y_k\in\mathcal{Y}$.

For arbitrary $k\ge k_{\epsilon}$, we have the following from (\ref{eqn:yk1}),(\ref{eqn:y_inf}).
\begin{equation*}
\begin{aligned}
&\|y_k-y_\infty\|_\lambda\\
=&\|G_k(y_{k-1})-G_\infty(y_\infty)+T_k(x_{k-1}-z_{k-1})\|_\lambda\\
\le & \|G_k(y_{k-1})-G_k(y_\infty)\|_\lambda+\|G_k(y_\infty)-G_\infty(y_\infty)\|_\lambda\\
&+\|T_k(x_{k-1}-z_{k-1})-T_k(0)\|_\lambda\\
\le & \rho\|y_{k-1}-y_\infty\|_\lambda+2\epsilon'D_2+\rho D_1
\end{aligned}
\end{equation*}
Iterating the above equation to $k_{\epsilon}$, we have
\begin{equation*}
\begin{aligned}
\|y_k-y_\infty\|_\lambda\le &\rho^{k-k_{\epsilon}}\|y_{k_{\epsilon}}-y_\infty\|_\lambda+\sum_{i=0}^{{k-k_{\epsilon}-1}}2\rho^i\epsilon'D_2\\
\le&\rho^{k-k_{\epsilon}}\|y_{k_{\epsilon}}-y_\infty\|_\lambda+\frac{2\epsilon'D_2}{1-\rho}+\frac{\rho D_1}{1-\rho}
\end{aligned}
\end{equation*}
It follows from Theorem \ref{thm} that $\|y_{k\epsilon}-y_\infty\|_\lambda$ is bounded by $2D_2$. Therefore, it can be guaranteed that $\|y_k-y_\infty\|_\lambda\le D_y+\epsilon$ for any $k$ satisfying
\begin{equation*}
k>k_{\epsilon}+\log_{\rho}\frac{1}{2D_2}\left(\epsilon-\frac{2\epsilon'D_2}{1-\rho}\right)
\end{equation*}
$\epsilon$ is chosen arbitrarily; thus, we can conclude that
\begin{equation*}
\lim_{k\to\infty}\|y_k-y_\infty\|_\lambda\le \frac{\rho D_1}{1-\rho}
\end{equation*}
It follows that
\begin{equation*}
\lim_{k\to\infty}\text{dist}(y_k,\mathcal{Y})= 0
\end{equation*}
That means $y_k$ will finally be within $\mathcal{Y}$. This completes the proof.
\end{IEEEproof}

%\section{}

\section*{Acknowledgment}

The authors would like to thank Prof. John Hunter at the University of California, Davis for the helpful discussion, and the associate editor and the reviewers for a number of helpful comments based on which the paper has been greatly enhanced. %Academician Workstation Project of Guangdong, China $\#$2012B090500010, National Natural Science Foundation of China $\#$61433005 and Hong Kong Research Grant Council Project $\#$612512 for supporting this work. 

\ifCLASSOPTIONcaptionsoff
  \newpage
\fi



%\bibliographystyle{IEEEtran}
%\begin{thebibliography}{1}
%\end{thebibliography}
\bibliographystyle{IEEEtran}
\bibliography{autosam}
%\begin{IEEEbiography}[{\includegraphics[width=1in,height=1.25in,clip,keepaspectratio]{Edward}}]{Zhixing Cao (Edward)}
%\end{IEEEbiography}

%\begin{IEEEbiography}[{\includegraphics[width=1in,height=1.25in,clip,keepaspectratio]{Edward}}]{Zhixing Cao (Edward)}
\begin{IEEEbiographynophoto}{Zhixing Cao (Edward)}
received his B.Eng. from the Department of Control Science and Engineering, Zhejiang University, China, in 2012 and his Ph.D. degree in Chemical and Biomolecular Engineering, from the Hong Kong University of Science and Technology (HKUST) in 2016. He is now a postdoctoral fellow in Harvard John A. Paulson School of Engineering and Applied Sciences, Harvard University. He is also a recipient of the Hong Kong PhD Fellowship. His research interests include iterative learning control, system identification, extremum seeking control and their applications in batch process and bioengineering.
\end{IEEEbiographynophoto}

%\begin{IEEEbiography}[{\includegraphics[width=1in,height=1.25in,clip,keepaspectratio]{durr}}]{Hans-Bernd D{\"u}rr} 
\begin{IEEEbiographynophoto}{Hans-Bernd D{\"u}rr}
studied at the University of Stuttgart, Germany, at the {\'E}cole Centrale Paris, France and at the Royal Institute of Technology (KTH), Sweden. From the University of  Stuttgart he received a Diploma degree in 2010 and a PhD in 2015 in the field of Engineering Cybernetics. He is currently working as a development engineer in automotive industry. His research interests include optimization and control theory.
\end{IEEEbiographynophoto}

%\begin{IEEEbiography}[{\includegraphics[width=1in,height=1.25in,clip,keepaspectratio]{ebenbauer}}]
\begin{IEEEbiographynophoto}{Christian Ebenbauer} received his MS (Dipl.-Ing.) in Telematics (Electrical Engineering and Computer Science) from Graz University of Technology, Austria, in 2000 and his Ph.D. (Dr.-Ing.) in Mechanical Engineering from the University of Stuttgart, Germany, in 2005. After having completed his Ph.D., he was a Postdoctoral Associate and an Erwin Schr{\"o}dinger Fellow at the Laboratory for Information and Decision Systems, Massachusetts Institute of Technology, USA. Since April 2009, he is a full professor at the Institute for Systems Theory and Automatic Control, University of Stuttgart, Germany. His research interests lie in the areas of dynamical systems, control theory, optimization and computation.
\end{IEEEbiographynophoto}

%\begin{IEEEbiography}[{\includegraphics[width=1in,height=1.25in,clip,keepaspectratio]{allgower}}]
\begin{IEEEbiographynophoto}{Frank Allg{\"o}wer} is the director of the Institute for Sys- tems Theory and Automatic Control at the University of Stuttgart. He studied Engineering Cybernetics and Applied Mathematics in Stuttgart and at UCLA respectively and received his Ph.D. from the University of Stuttgart. Prior to his present appointment he held a professorship in Electrical Engineering at ETH Zurich. He received several recognitions for his work including the prestigious Gottfried-Wilhelm-Leibniz prize of the Deutsche Forschungsgemeinschaft. His main areas of interest are in cooperative control, predictive control and systems biology.
\end{IEEEbiographynophoto}

%\begin{IEEEbiography}[{\includegraphics[width=1in,height=1.25in,clip,keepaspectratio]{Furong_Gao}}]
\begin{IEEEbiographynophoto}{Furong Gao}is a Chair Professor in the Department of Chemical and Biomolecular Engineering at the Hong Kong University of Science and Technology (HKUST). He obtained his B.Eng. in Automation from East China Institute of Petroleum in 1985, and his M.Eng. and Ph.D. degrees in Chemical Engi- neering from McGill University, Montreal, Canada, in 1989 and 1993, respectively. He worked as a Senior Research Engineer at Moldflow International, Melbourne, Australia, from 1993 to 1995 before joining HKUST as a professor. His research interests
include process monitoring and fault diagnosis, batch process control, polymer processing control, and optimization. He received a number of best paper awards, and is on Editorial Boards of a number of numerous journals of his area.
\end{IEEEbiographynophoto}
%% insert where needed to balance the two columns on the last page with
%%\newpage


%\vfill


\end{document}